\magnification=\magstephalf
\input eplain
\indexproofingtrue
\input epsf
\input BM.mac


\centerline{{\titlefont RELATIONS BETWEEN TAUTOLOGICAL CYCLES}}
\Bskip

\centerline{{\titlefont ON JACOBIANS}}
\Askip

\centerline{{\it by}}
\Bskip

\centerline{{\namefont Ben Moonen}}
\Askip\Askip

{\narrower {\eightpoint\noindent{\bf Abstract.} We study tautological cycle classes on the Jacobian of a curve. We prove a new result about the ring of tautological classes on a general curve that allows, among other things, easy dimension calculations and leads to some general results about the structure of this ring. Next we obtain a vanishing result for some of the generating classes~$p_i$; this gives an improvement of an earlier result of Herbaut. Finally we lift a result of Herbaut and van der Geer-Kouvidakis to the Chow ring (as opposed to its quotient modulo algebraic equivalence) and we give a method to obtain further explicit cycle relations. As an ingredient for this we prove a theorem about how Polishchuk's operator~$\cD$ lifts to the tautological subalgebra of~$\Chow(J)$.
\smallskip

\par}}
\Bskip


\introduction

\introsec{introbegin}
Let $(J,\theta)$ be the Jacobian of a curve~$C$ of genus $g \geq 2$, and let $j \colon C \hookrightarrow J$ be the embedding obtained by choosing a point $x_0 \in C$. The Chow ring $\Chow(J)$ (with $\mQ$-coefficients) comes equipped with a number of structures: in addition to the usual intersection product we have the Pontryagin product~$*$, the Fourier transform~$\cF$, and the action of~$n^*$ and~$n_*$ for all $n \in \mZ$. Using the latter we obtain Beauville's decomposition $\Chow(J) = \oplus\, \Chow^i_{(j)}(J)$. These structures are inherited by the quotient $A(J) := \Chow(J)/ \sim_\alg$.

In [\refn{BeauvACJ}], Beauville studied the tautological subring $\cT(C) \subset A(J)$, that can be defined as the smallest subalgebra that is stable under all operations just mentioned, and that contains the image of $A(C)$ under~$j_*$. He proved that $\cT(C)$ is generated by certain classes~$p_m$ that are the components of $\cF[C]$ in the Beauville decomposition. This leads to the question what is the ideal of relations between these classes.

The work of Polishchuk~[\refn{PolUniv}] provides us with a powerful method to produce relations. He considers the map $\phi\colon R \twoheadrightarrow \cT(C)$ from the polynomial ring $R := \mQ[x_1,x_2,\ldots]$ given by $x_i \mapsto p_i$. The elements of~$R$ that are of degree~$\geq g$ (with $x_i$ of degree~$i$) with zero coefficient in front of~$x_1^g$ lie in the kernel of~$\phi$ for obvious reasons. But now the point is that Polishchuk is able to give an explicit differential operator~$\cD$ on~$R$ that preserves the kernel of~$\phi$, and that the ideal of ``trivial'' relations is far from stable under~$\cD$. Polishchuk studies the smallest ideal $I \subset R$ that is $\cD$-stable and that contains all trivial relations. Conjecturally, for general curves~$C$ (over a base field of sufficiently high transcendence degree) the quotient ring $\cR := R/I$ maps isomorphically to the tautological ring~$\cT(C)$.

The present paper has three main goals:
\Bskip

\noindent
(1) First we want to give some new results on the structure of the ring~$\cR$. Though in principle it is not difficult to calculate, for a given~$g$, the dimensions of all summands~$\cR^i_{(j)}$ on a computer, it is not so easy to obtain general conclusions about the structure of~$\cR$ from Polishchuk's methods. (In some sense the main obstacle is that the variable~$x_1$ plays an exceptional role in the calculations.)

We prove a theorem, Thm.~\refn{RjStructThm}, that gives more insight in the structure of~$\cR$. Each $\cR_{(j)} := \oplus_i\, \cR^i_{(j)}$ is a finite-dimensional $\sltwo$-module. We give a simple recipe for the multiplicity of a given irreducible $\sltwo$-module in each~$\cR_{(j)}$. As corollaries we get some general results about the relations in~$\cR_{(j)}$, especially for small and large (relative to~$g$) values of~$j$. Further we present a conjecture of van der Geer and Kouvidakis that gives a simple recipe for the dimensions of the spaces~$\cR^i_{(j)}$. 
\Bskip

\noindent
(2) The second goal of the paper is to investigate cycle relations for curves with a given linear system. If the curve has a~$g^r_d$, Herbaut~[\refn{Herbaut}] proved that we get new relations between the classes~$p_m$. This result was reproved and simplified by van der Geer and Kouvidakis in~[\refn{vdGK}]. However, the ideal generated by~$I$ together with these new relations is not stable under the operator~$\cD$. Ideally we should like to get a full description of the quotient of~$\cR$ by the ideal generated by Herbaut's relations together with all their images under powers of~$\cD$. This turns out to be not so easy. However, we do obtain a vanishing result, Thm.~\refn{DrminHerb} for classes~$p_k$ with $k > d+1-2r$. It turns out that the result we obtain is a slight improvement of an earlier result of Herbaut in~[\refn{Herbaut}].    
\Bskip

\noindent
(3) The third goal of the paper is to lift some results obtained thus far to the Chow ring of~$J$. For general curves Polishchuk already obtained a number of results about this in~[\refn{Polish}]. He considered a subalgebra $\taut(C) \subset \Chow(J)$ that we call the small tautological ring, and he proved that is generated by classes~$p_m$ (lifting those considered above) and~$q_m$ (essentially the Beauville components of~$\cF(j_* K)$, with $K$ the canonical class on~$C$).

In itself, it is not very hard to lift the result of Herbaut and van der Geer-Kouvidakis to the Chow level. We do this in Section~\refn{CRModAlg2}, following the method of~[\refn{vdGK}] which is based on a Grothendieck-Riemann-Roch calculation. See Thm.~\refn{ThmModRat} for the result.

However, if we also want to lift the operator~$\cD$ to the Chow level, we find that the ``small'' tautological ring needs to be enlarged, as the calculations involve natural classes that are in general not contained in it. So we define the big tautological ring $\Taut(C) \subset \Chow(J)$ as the smallest $\mQ$-subalgebra that is stable under all operations $\cdot$, $*$, $\cF$, $n^*$ and~$n_*$ and that contains the image of $j_* \colon \Chow(C) \to \Chow(J)$. Again there is an operator~$\cD$, but as $\Taut(C)$ is not finitely generated it is now less easy to describe this operator explicitely. However, the calculations we should like to do involve only finitely many classes of the form $\cF(j_* D_i)$ with $D_i$ divisors on~$C$. We show that the subalgebra of~$\Taut(C)$ generated by these classes together with the $p_m$ and~$q_m$ is stable under~$\cD$ and that $\cD$ again acts on it as a differential operator that we can give explicitely. Closely related results have been obtained independently by Alexander Polishchuk; see Remark~\refn{recentwork} below.

The situation we end up with is that we have lifted Herbaut's relations to the Chow level, and that we also have an explicit operator~$\cD$. Combined this gives an abundance of relations that can easily be calculated on a computer. We give some examples of the relations thus obtained in the final section. What seems much harder is to get some control over the total ideal of relations thus obtained. In particular, at this moment we have no analogue of Thm.~\refn{RjStructThm} for curves with a given linear system.  

\introsec{acknowledge}
{\bf Acknowledgement.} I thank Tom Koornwinder for his help getting started with Mathematica. I thank Gerard van der Geer and Alexis Kouvidakis for allowing me to present their conjecture on the dimensions of the spaces~$\cR^i_{(j)}$, see~\refn{ClosedForm}, and for explaining to me some details pertaining to their conjecture.

\introsec{recentwork}
{\bf Remark.} Very recently, some related results have appeared. In particular we should like to draw the reader's attention to Polishchuk's preprint~[\refn{PolRelSymm}], in which some ideas appear that are closely related to our Section~\refn{BigTaut}, and to the preprint~[\refn{FuHerbaut}] of Fu and Herbaut, that contains some results related to the material in Section~\refn{CRinChow}. 

\introsec{NotConv}
{\bf Notation and conventions.} We work over a fixed algebraically closed field~$k$. Throughout, if $X$ is a non-singular complete variety then by $\Chow(X)$ we denote the Chow ring of~$X$ {\it tensored with~$\mQ$.\/}

We write $\sltwo = \mQ \cdot f + \mQ \cdot h + \mQ \cdot e$ with $[e,f] = h$, $[e,h] = -2e$, and $[f,h] = 2f$.


\section{The Chow ring of a Jacobian}{ChowJac}

\ssection{situation}
Let $C$ be a non-singular complete curve of genus $g \geq 2$. Let $J := \Pic^0_{C/k}$ be its Jacobian, $\theta \colon J \isomarrow J^t$ the canonical polarization. We also write $\theta \in \Chow^1(J)$ for the class of a symmetric theta-divisor. Choose a base point $x_0 \in C$, and let $j \colon C \hookrightarrow J$ be the embedding given on points by $x \mapsto \bigl[O_C(x-x_0)\bigr]$.

Let $\cP_J$ be the Poincar\'e bundle on $J \times J^t$. Let $\cL_J$ be the line bundle on $J \times J$ obtained by pulling back~$\cP_J^{-1}$ under the map $\id_J \times \theta$. We have a Fourier transform $\cF \colon \Chow(J) \to \Chow(J)$ given by $\cF(x) = \pr_{2,*}\bigl(\chern(\cL_J) \cdot \pr_1^*(x)\bigr)$. Further we have Beauville's decomposition $\Chow^i(J) = \oplus_j\, \Chow^i_{(j)}(J)$, with $x \in \Chow^i_{(j)}(J)$ if and only if $n^*(x) = n^{2i-j}\cdot x$ for all $n \in \mZ$. The Fourier transform~$\cF$ induces a bijection between $\Chow^i_{(j)}(J)$ and $\Chow^{g-i+j}_{(j)}(J)$.

We write $\bigl[j(C)\bigr]$, or if there is no risk of confusion simply~$[C]$, for the class of $j(C)$ in $\Chow^{g-1}(J)$. Further, $[j(C)]_{(n)}$, or simply~$[C]_{(n)}$, denotes the component of~$[C]$ in $\Chow^{g-1}_{(n)}(J)$.

\ssection{picture}
We find it helpful to draw pictures representing the Chow ring as in Figure~1.
\topinsert 
$$
\epsfxsize=\hsize\epsfbox{chowpicture.1}
$$
\centerline{{\bf Figure~1. A picture of $\Chow(J)$.}}
\endinsert
The boxes $(i,j)$ represent the spaces $\Chow^i_{(j)}$. In the horizontal direction we have the {\it weight\/}, where by definition $\Chow^i_{(j)}$ has weight $2i-j$. The Fourier transform acts as reflection in the central vertical line.

The filtration $\Fil^\gdot$ on $\Chow(J)$ defined by $\Fil^r\Chow := \oplus_{j \geq r}\, \Chow^\gdot_{(j)}$ should satisfy the conjectures of Beilinson and Murre; see Jannsen~[\refn{JannMotSh}]. Most of the expected properties of this filtration are as yet unproved; we shall mention some of these.

It is known that $\Chow^i_{(j)}$ can be nonzero only for $g-i \leq j \leq i$. According to a conjecture of Beauville in [\refn{BeauvQuelq}], \S~5, the spaces $\Chow^i_{(j)}$ should be zero if $j < 0$; this corresponds to one of the properties of Beilinson's conjectural filtration. In general this is known only for $i \in \{0,1,g-2,g-1,g\}$. Over finite fields or the algebraic closure of a finite field Beauville's conjecture is known; see K\"unnemann [\refn{Kunne}], Sections 7 and~8. In the rest of the paper we shall usually omit the ``negative~$j$'' part in our drawings.

As $n^*$ acts on $H^{2i}$ as multiplication by~$n^{2i}$, all classes in $\Fil^1\Chow$ are homologically trivial. It should be the case that $\oplus_i\, \Chow^i_{(0)}$ injects into cohomology, but this is not known in general. Similarly, by considering the weight we see that all classes in $\Fil^2\Chow$ map to zero under the Abel-Jacobi map, and it should be the case that $\Fil^2\Chow$ is precisely the kernel of the Abel-Jacobi map.

The classes in $\Fil^1\Chow^g$ represent $0$-cycles of degree zero and are therefore algebraically trivial. By Fourier duality it follows that also all classes in $\Chow^j_{(j)}$ with $j > 0$ are algebraically trivial. In the picture these summands are indicated by boxes with a heavier border. As we shall see, in general there are many more classes that are algebraically trivial.

\ssection{sl2onchow}
By the work of K\"unnemann in~[\refn{Kunne}], the Chow motive of~$J$ has a Lefschetz decomposition. This gives rise to an action of $\sltwo$ on~$\Chow(J)$. We normalise this as in [\refn{Polish}]; so with notation as in~\refn{NotConv} we have
$$
e(\alpha) = -\theta \cdot \alpha\, ,
\qquad
f(\alpha) = [C]_{(0)} * \alpha\, ,
$$ 
and
$$
h(\alpha) = (2i-j-g) \cdot \alpha\qquad \hbox{for $\alpha \in \Chow^i_{(j)}(J)$.}
$$

\ssection{chowmodalg}
Define $A(J)$ as the quotient of $\Chow(J)$ modulo algebraic equivalence. This ring inherits all the structures on $\Chow(J)$ that are relevant for us. Concretely, $A(J)$ has an intersection product, a Pontryagin product, a Fourier transform, a Beauville decomposition $A(J) = \oplus\, A^i_{(j)}$, and an $\sltwo$-action. We shall use these structures without further comments. We again write $\theta$ for the class of a symmetric theta-divisor in $A^1_{(0)}(J)$.


\section{Cycle relations modulo algebraic equivalence---general curves}{CRmodAlg1}

\ssection{theringR}
Throughout this section, $g$ will be a fixed integer with $g \geq 3$. We start by reviewing some results from Polishchuk's paper~[\refn{PolUniv}]. 

Consider the polynomial ring $R := \mQ[x_1,x_2,\ldots]$ in infinitely many variables. There are several gradings that will play a role in what follows. Among these are gradings that we call ``codimension'', ``level'' and ``weight''; they are defined by setting
$$
\codim(x_i) := i\, ,\qquad
\ell(\alpha) := i-1\, ,
\qquad\hbox{and}\qquad
w(x_i) := i+1\, .
$$
We have $w = 2\cdot \codim - \ell$.

Define $R^i_{(j)} := \bigl\{\alpha \in R \bigm| \hbox{$\codim(\alpha) = i$ and $\ell(\alpha) = j$} \bigr\}$. This gives a bigrading $R = \oplus\, R^i_{(j)}$. We write $R^i := \oplus_j\, R^i_{(j)}$ and $R_{(j)} := \oplus_i\, R^i_{(j)}$. Let $\Fil^\gdot$ be the descending filtration by level, so $\Fil^n R := \oplus_{j \geq n}\, R_{(j)}$.

Define an operator $\cD$ on~$R$ by
$$
\cD := -g\partial_1 + {\textstyle{1\over 2}} \sum_{m,n \geq 1}\, {m+n\choose n} x_{m+n-1} \partial_m\partial_n\, , 
$$
where $\partial_k := \partial_{x_k}$. We give $R$ the structure of an $\sltwo$-module by setting
$$
\eqalign{
e(\alpha) &= x_1 \cdot \alpha\, ;\cr
h(\alpha) &= \bigl(w(\alpha)-g\bigr) \cdot \alpha = (2i-j-g) \cdot \alpha \qquad\hbox{if $\alpha \in R^i_{(j)}$}\, ;\cr
f(\alpha) &= -\cD(\alpha)\, .}
$$
See [\refn{PolUniv}], Lemma 3.2, and see~\refn{NotConv} for our notation regarding~$\sltwo$.

Define $I \subset R$ to be the smallest linear subspace that is stable under~$\cD$ and that contains $R^{>g} + \Fil^1 R^g$. Concretely, $I = R^{>g} \oplus \sum_{n \geq 0}\, \cD^n\bigl(\Fil^1 R^g \bigr)$. Polishchuk shows that $I = \cap_{n \geq 0}\, \Image(\cD^n)$, and that $I$ is in fact an ideal of~$R$. Define $\cR := R/I$. (This is the ring called $R^\Jac_g$ in~[\refn{PolUniv}], but for later use it will be convenient to have a simpler notation.) The ring $\cR$ inherits a bigrading $\cR = \oplus\, \cR^i_{(j)}$ and the structure of an $\sltwo$-module. Note that the subspaces $\cR_{(j)}$ are $\sltwo$-submodules of~$\cR$.

The next ingredient is that we have a ``Fourier operator''~$\cF$ on~$\cR$, given by [\refn{PolUniv}], formula~(0.4). (Our $\cF$ is Polishchuk's~$S$.) If $\alpha \in \cR^i_{(j)}$ then $\cF^2(\alpha) = (-1)^{g+j} \alpha$.

\ssection{RtoChow}
Let now $C$ be a non-singular curve of genus~$g$. We use the notation of~\refn{situation}. Consider the homomorphism of $\mQ$-algebras $\phi \colon R \to A(J)$ given by $\phi(x_i) = p_i$, with
$$
p_i := \hbox{degree $n$ component of $\cF[C]$} = \cF\bigl([C]_{(n-1)}\bigr)\, .
$$
The image of~$\phi$ is the tautological subring introduced and studied by Beauville in~[\refn{BeauvACJ}]. We shall denote this tautological ring by $\cT(C) := \Image(\phi) \subset A(J)$.

We have $I \subset \Ker(\phi)$. Let $\psi \colon \cR \to A(J)$ be the induced homomorphism. This map $\psi$ is compatible with all structures considered above; that is, \item{---} $\psi(\cR^i_{(j)}) \subseteq A^i_{(j)}(J)$ for all $i$, $j$;
\item{---} $\psi$ is $\sltwo$-equivariant;
\item{---} $\psi \circ \cF = \cF \circ \psi$.

\ssection{OurGoal}
The tautological ring $\cT(C) \subset A(J)$ is generated by the classes $p_1,\ldots, p_g$. We should like to understand the relations between these classes. The relations in Polishchuk's ideal~$I$ are those that are obtained from ``trivial relations'', i.e., relations $F(p_1,p_2,\ldots) = 0$ with $F \in R^{>g} + \Fil^1 R^g$, by (repeated) application of the operator~$\cD$. The relations thus obtained are in many cases non-trivial from a geometric perspective. For instance, using this method Polishchuk shows that $p_n = 0$ in~$A(J)$ for all $n \geq {g\over 2} + 1$. (These particular relations also follow from a theorem of Colombo and van Geemen in~[\refn{ColvGeem}], as the gonality of~$C$ is at most $1 + \lceil{g/2}\rceil$.) Polishchuk conjectures that for a generic curve~$C$ (over an algebraically closed base field of sufficiently high transcendence degree over the prime field) the ideal~$I$ gives all relations.

For a given~$g$ we can, at least in principle, write down a basis for the ideal~$I$ over~$\mQ$. See [\refn{PolUniv}], section~2.6 for examples in genera $\leq 10$. The drawback of this method is that it is purely computational, and that it is hard to get an insight in the structure of~$\cR$. For instance, in examples one finds that there are many pairs $(i,j)$, especially $j \geq g/2$, for which $\cR^i_{(j)} = 0$, and one would like to get an insight in when this occurs.

Our main result in this section further details Polishchuk's method, and gives the structure of the ring~$\cR$ as an $\sltwo$-module. Using our theorem, it becomes very easy to calculate the dimensions the~$\cR^i_{(j)}$, and we also get some general results about the structure of~$\cR$.

We start with a simple lemma.

\ssection{Dnp1aalpha}
{\bf Lemma.}\enspace {\it Let $\alpha \in R^i_{(j)}$. Then 
$$
\cD^n(x_1^a \cdot \alpha) = \sum_{s=0}^{\min(n,a)}\; {n!\over (n-s)!} {a! \over (a-s)!} \cdot {2i-j-g+a-n+s-1\choose s} \cdot\; x_1^{a-s} \, \cD^{n-s}(\alpha)
$$
for all $n \geq 0$ and $a \geq 0$.}
\Cskip

Note that the binomial coefficient has to be taken in the generalised sense, as $2i-j-g+a-n+s-1$ may be negative. As a particular case of the lemma, we have
$$
\cD(x_1^a \cdot \alpha) = x_1^a \cdot \cD(\alpha) + a(2i-j-g+a-1) \cdot x_1^{a-1} \alpha\, .\eqlabel{Dp1aalpha}
$$
\Dskip

\Proof~The formula is proven by an easy induction on~$n$. To start the induction we first prove (\refn{Dp1aalpha}) by induction on~$a$, using that $\cD(x_1 \beta) - x_1\cD(\beta) = (2i-j-g) \cdot \beta$ for $\beta \in R^i_{(j)}$. \QED

\ssection{Monomials}
Given an integer $j \geq 0$, define $\Mon_{(j)} \subset R_{(j)}$ to be the (finite) set of monomials in the variables $x_2, x_3, \ldots$ (no~$x_1$) with $\ell(\alpha) = j$. For $j=0$ we have $\Mon_{(j)} = \{1\}$. Let $\Mon^i_{(j)} := \Mon_{(j)} \cap R^i_{(j)}$.

If $j \geq 1$ there is a bijection between $\Mon_{(j)}$ and the set of ordered partitions $j = j_1 + j_2 + \cdots + j_r$ with $1 \leq j_1 \leq j_2 \leq \cdots \leq j_r$, letting such a partition correspond to $\alpha = x_{j_1+1} x_{j_2+1} \cdots x_{j_r+1}$. Under this bijection, $\Mon^i_{(j)}$ corresponds with the partitions with $i-j$ parts. In particular, $\Mon^i_{(j)}$ is non-empty (still for $j\geq 1$) if and only if $j+1 \leq i \leq 2j$. 

The set $\Mon_{(j)}$ is a basis for $R_{(j)}$ as a $\mQ[x_1]$-module. Write $M^i_{(j)} := \mQ \cdot \Mon^i_{(j)}$ and $M_{(j)} := \mQ \cdot \Mon_{(j)}$. Note that $\cD(M^i_{(j)}) \subseteq M^{i-1}_{(j)}$.

\ssection{RjStructThm}
{\bf Theorem.}\enspace {\it Write $\Standard$ for the tautological $2$-dimensional representation of~$\sltwo$. Given $i$ and~$j$ with $2i-j \leq g$, write
$$
\mu(i,j) := \dim_\mQ\bigl(M^i_{(j)}/\cD^{g-2i+j+1}(M^{g-i+j+1}_{(j)})\bigr)\, .
$$
Then for all $j\geq 0$ we have
$$
\cR_{(j)} \cong \bigoplus_{i=j}^{\min(g-1,2j,\lfloor{g+j\over 2}\rfloor)}\; \bigl[ \Sym^{g-2i+j}(\Standard)\bigr]^{\mu(i,j)}
\eqlabel{cRStructIsom}
$$
as $\sltwo$-modules.}
\Dskip

\Proof~The result for $j=0$ says that $\cR_{(0)} \cong \Sym^g(\Standard)$, which is immediate. In the rest of the proof we shall assume that $j \geq 1$. Let $\Phi_\gdot$ be the ascending filtration of $R_{(j)}$ by $\mQ[x_1]$-submodules that is given by $\Phi_i := \mQ[x_1] \cdot \bigl(\oplus_{k \leq i}\,  M^k_{(j)}\bigr)$. Note that $\Phi_j = (0)$ and $\Phi_{2j} = R_{(j)}$. Let $\Psi_\gdot$ be the filtration on~$\cR_{(j)}$ induced by~$\Phi_\gdot$.

If $\alpha \in \Mon^i_{(j)}$ then $\cD^n(\alpha) \in \Phi_{i-1}$ for all $n \geq 1$. Hence by Lemma~\refn{Dnp1aalpha} we have
$$
\cD^n(x_1^a \cdot \alpha) \equiv \cases{0 & if $a < n$;\cr
{a!\over (a-n)!} {2i-j-g+a-1\choose n}\; x_1^{a-n} \alpha & if $a \geq n$;\cr}
$$
modulo $\Phi_{i-1}$. Write $V_\alpha$ for the image of $\mQ[x_1] \cdot \alpha$ in $\gr_i^\Phi$, which is an $\sltwo$-submodule. We find:
\item{---} If $2i-j > g$ then $V_\alpha$ is irreducible, $\dim(V_\alpha) = \infty$.
\item{---} If $2i-j \leq g$ then $V_\alpha$ is an extension,
$$
0 \tto U_\alpha \tto V_\alpha \tto W_\alpha \tto 0\, ;
$$
here $U_\alpha := \Image\bigl(\mQ[x_1] \cdot x_1^{g-2i+j+1} \alpha\bigr)$ is infinite dimensional and irreducible, and $W_\alpha$ is isomorphic with $\Sym^{g-2i+j}(\Standard)$. 
\smallskip

\noindent
As $\cR_{(j)}$ is finite dimensional it follows in particular that $\Phi_{\lfloor{g+j\over 2} \rfloor}$ surjects to~$\cR_{(j)}$.

Let $i$ be an integer with $j+1 \leq i \leq \min\bigl({2j,\lfloor{g+j\over 2}\rfloor}\bigr)$. The set $\Mon^i_{(j)}$ gives a basis for $\gr_i^\Phi$ as a $\mQ[x_1]$-module, and by what we have just seen we have an exact sequence
$$
0 \tto U \tto \gr_i^\Phi \tto W \tto 0\, ,
$$
where $W$ is the direct sum of the spaces $W_\alpha$ for $\alpha \in \Mon^i_{(j)}$. Consider the subspace $M^i_{(j)} \subset \Phi_i$. Let $\bar{M}^i_{(j)}$ be the image of~$M^i_{(j)}$ under the composition $\Phi_i \twoheadrightarrow \gr_i^\Phi \twoheadrightarrow W$. The natural map $M^i_{(j)} \to \bar{M}^i_{(j)}$ is an isomorphism. With $f$ as in~\refn{NotConv} we have $\bar{M}^i_{(j)} = \Ker(f_{|W})$, and $W \cong \bar{M}^i_{(j)} \otimes_\mQ \Sym^{g-2i+j}(\Standard)$ as $\sltwo$-modules.

Clearly the natural map $\gr_i^\Phi \twoheadrightarrow \gr_i^\Psi$ factors via a map $\xi \colon W \twoheadrightarrow \gr_i^\Psi$. Let $K := \Ker(\xi)$, which is an $\sltwo$-submodule of~$W$. Let $K[f] := \Ker(f_{|K})$. Then we have $K[f] = \Ker(\xi) \cap \bar{M}^i_{(j)}$, and $\gr_i^\Psi$ is isomorphic with $\bigl(\bar{M}^i_{(j)}/K[f]\bigr) \otimes_\mQ \Sym^{g-2i+j}(\Standard)$ as $\sltwo$-modules.

Now consider the composition
$$
M^i_{(j)} \isomarrow \bar{M}^i_{(j)} \longhookrightarrow W \twoheadrightarrow \gr_i^\Psi\, ,
$$
which is the just the restriction of the natural map $\Phi_i \twoheadrightarrow \gr_i^\Psi$ to $M^i_{(j)} \subset \Phi_i$. An element $y \in M^i_{(j)}$ maps to zero under this map if and only if $y \in I + \Phi_{i-1}$. So to complete the proof, it suffices to show that for $y \in M^i_{(j)}$ we have
$$
\eqalignno{
y \in I + \Phi_{i-1} &\quad\Leftrightarrow\quad  y \in \bigl(I \cap R^i_{(j)}\bigr) + x_1 R^{i-1}_{(j)} & {\rm (a)}\cr
&\quad\Leftrightarrow\quad  y \in \cD^{g-i}(R^g_{(j)}) + x_1 R^{i-1}_{(j)} & {\rm (b)}\cr
&\quad\Leftrightarrow\quad  y \in \left(\sum_{a=0}^{i-j-1} \cD^{g-i-a}(M^{g-a}_{(j)})\right) + x_1 R^{i-1}_{(j)} & {\rm (c)}\cr
&\quad\Leftrightarrow\quad  y \in \cD^{g-2i+j+1}(M^{g-i+j+1}_{(j)}) + x_1 R^{i-1}_{(j)} & {\rm (d)}\cr
&\quad\Leftrightarrow\quad  y \in \cD^{g-2i+j+1}(M^{g-i+j+1}_{(j)}) \, .& {\rm (e)}\cr}
$$ 

For (a), the implication ``$\Leftarrow$'' is clear, as $x_1 R^{i-1}_{(j)} \subseteq \Phi_{i-1}$ and $I \cap R^i_{(j)} \subseteq I$. For the converse, suppose we have $z \in I$ and $w \in \Phi_{i-1}$ with $y = z+w$. Because $I$ and $\Phi_{i-1}$ are bi-homogeneous with respect to the decomposition $R = \oplus R^i_{(j)}$, we may replace $z$ and~$w$ by their components in~$R^i_{(j)}$, in which case $z \in I \cap R^i_{(j)}$ and $w \in x_1 \cdot R^{i-1}_{(j)}$.

For (b), just note that $I \cap R^i_{(j)} = \cD^{g-i}(R^g_{(j)})$. For~(c), which is really the main point, start with the decomposition $R^g_{(j)} = \oplus_{a \geq 0}\, x_1^a \cdot M^{g-a}_{(j)}$. If $\alpha \in M^{g-a}_{(j)}$ then using Lemma~\refn{Dnp1aalpha} we find that
$$
\cD^{g-i}(x_1^a \alpha) = \cases{0 & if $a \geq i-j$;\cr
{(g-i)!\; a!\over (g-i-a)!} \; {i-j-1\choose a}\; \cD^{g-i-a}(\alpha) & for $0 \leq a \leq i-j-1$,}
$$
modulo $x_1 \cdot R^{i-1}_{(j)}$.

Finally, (d) follows from the remark that $\cD(M^b_{(j)}) \subseteq M^{b-1}_{(j)}$ for all~$b$, and for (e) we use that $y$ and $\cD^{g-2i+j+1}(M^{g-i+j+1}_{(j)})$ are both contained in $M^i_{(j)}$, whereas $M^i_{(j)} \cap  x_1 R^{i-1}_{(j)} = 0$. \QED

\ssection{RjStructRem}
{\bf Remarks.\/}\enspace (\romno1) One of the main advantages of our result is that the calculations do not involve the variable~$x_1$, and that therefore the operator~$\cD$ becomes a lot easier. In fact, we feel the corollaries below do not yet represent the best possible conclusions that one should be able to get. See \refn{ClosedForm} for some speculation. The only obstacle for pushing our results further is of a purely combinatorical nature.

(\romno2) The proposition should not be read as saying that the images of the monomial spaces $M^i_{(j)}$ in $\cR^i_{(j)}$ consist of primitive classes, i.e., classes in the kernel of~$f$. This is simply not the case. However, to an element $\alpha \in M^i_{(j)}$ corresponds a primitive class~$\alpha^\prime$, given by
$$
\alpha^\prime = \sum_{n\geq 0}\; {(-x_1)^n \cdot \cD^n(\alpha)\over (n!)^2\cdot {2i-j-g-2\choose n}}\, .
$$
(Note that the sum is finite, as $\cD^n(\alpha) = 0$ for $n$ large.)

(\romno3) For $j\geq 1$ the direct sum in~(\refn{cRStructIsom}) starts at $i=j+1$. To ``visualize'' where the various summands $M^i_{(j)}$ and $R^i_{(j)}$ lie, for instance in the definition of the multiplicity~$\mu(i,j)$, it is usually helpful to look at a picture, as in Figure~1. Recall that the number $2i-j$ is the weight. The summand $R^{g-i+j}_{(j)}$ is the Fourier mirror image of~$R^i_{(j)}$. The reader is encouraged to look at the examples in~\refn{lowgExa} below.

\ssection{ghigherg}
{\bf Corollary.} {\it Write $\cR[g]$ for the ring~$\cR$ in genus~$g$. If $2i-j \leq g$ and $g < g^\prime$ then the multiplicity of $\Sym^{g^\prime-2i+j}(\Standard)$ in~$\cR[g^\prime]_{(j)}$ is greater or equal to the multiplicity of $\Sym^{g-2i+j}(\Standard)$ in~$\cR[g]_{(j)}$.}
\Dskip

\Proof~Note that $\cD^{g^\prime-2i+j+1} \colon M^{g^\prime-i+j+1}_{(j)} \to M^i_{(j)}$ factors through $\cD^{g-2i+j+1} \colon M^{g-i+j+1}_{(j)} \to M^i_{(j)}$. \QED

\ssection{pnnonvanish}
{\bf Corollary.} {\it For $n \leq {g+1\over 2}$ we have $x_n \neq 0$ in~$\cR$.}
\Dskip

\Proof
We need to prove that for $j \leq {g-1\over 2}$ the representation $\Sym^{g-2-j}(\Standard)$ occurs in $\cR_{(j)}$. This corresponds to $i=j+1$ in~(\refn{cRStructIsom}). Now use that $M^g_{(j)}$ is zero for $j \leq {g-1\over 2}$. \QED
\Bskip

\noindent
Of course it is far more interesting to have results about the non-vanishing of classes~$p_n$ in $\Chow(J)$, for a general curve~$C$. Ceresa's theorem in~[\refn{Ceresa}] gives such a result, as it tells us that for $g \geq 3$ and general~$C$ we have $p_2 \neq 0$ in $A(J)$. (This follows from~[\refn{Ceresa}] using [\refn{Herbaut}], Thm.~5.) See Fakhruddin~[\refn{Fakh}] and Ikeda~[\refn{Ikeda}] for some further results.

\ssection{RjStructCor}
{\bf Corollary.\/}\enspace {\it {\rm (\romno1)} If $i+j \leq g$ then the natural map $R^i_{(j)} \to \cR^i_{(j)}$ is an isomorphism. 

{\rm (\romno2)} If $0 < j \leq g/3$ then $R^i_{(j)} \cong \cR^i_{(j)}$ for all~$i$ with $2i-j \leq g$. In this case we have
$$
\cR_{(j)} \cong \bigoplus_{i=j+1}^{2j}\; 
\bigl[\Sym^{g-2i-j}(\Standard)\bigr]^{\nu(i,j)} 
$$
as $\sltwo$-modules, where $\nu(i,j)$ is the number of ordered partitions of~$j$ with $i-j$ parts.}
\Dskip

\Proof~(\romno1) As remarked in~\refn{Monomials}, $M^b_{(j)} = 0$ if $b > 2j$. So if $i+j \leq g$ then in the theorem we have $\mu(i,j) = \dim(M^i_{(j)})$. This means that the kernel of the map $R_{(j)} \twoheadrightarrow \cR_{(j)}$ is contained in $R^{\geq g-j+1}_{(j)}$.

(\romno2) If $j \leq g/3$ then (\romno1) applies to all~$i$ in the decomposition~(\refn{cRStructIsom}). So then $\mu(i,j) = \dim_\mQ(M^i_{(j)})$, which, as we have seen in~\refn{Monomials}, is the number of ordered partitions of~$j$ with $i-j$ parts. \QED
\Bskip

\noindent
The previous corollaries can be interpreted as a non-vanishing results. The picture that emerges is that for small~$j$ there are no relations in~$R_{(j)}$. By contrast, for large~$j$ we expect many relations. We always have $\cR_{(j)} = 0$ for $j \geq g-1$, and Polishchuk's result about the vanishing of the class of~$x_j$ for $j \geq {g\over 2} +1$ (see [\refn{PolUniv}], Cor.~0.2) gives that $\cR_{(g-2)} = 0$ for all $g \geq 4$ and $\cR_{(g-3)} = 0$ for all $g \geq 6$. We expect that for a given~$j$ there is a bound~$G_j$ such that $\cR_{(j)} = $ for all $g \geq G_j$. In fact, the conjecture of van der Geer and Kouvidakis---see \refn{ClosedForm} below---predicts that $\cR_{(g-2l)} = 0$ for all $g \geq (l+1)^2$ and $\cR_{(g-2l+1)} = 0$ for all $g \geq l(l+1)$.  The following result gives a proof of this in the first non-trivial cases. The bounds we obtain are sharp; see the examples in~\refn{lowgExa}.

\ssection{j=gmin4}
{\bf Corollary.\/} {\it For all $g \geq 9$ we have $\cR_{(g-4)} = 0$. For all $g \geq 12$ we have $\cR_{(g-5)} = 0$.}
\Dskip

\Proof For the first assertion, assume $g \geq 9$. By [\refn{PolUniv}] we already know that $p_{g-3} = 0$; this means that $\Sym^3(\Standard)$ does not occur in $\cR_{(g-4)}$. Hence to prove that $\cR_{(g-4)} = 0$ it suffices, by the theorem, to show that the map $\cD \colon M_{(g-4)}^{g-1} \to M_{(g-4)}^{g-2}$ is surjective. A basis for $M_{(g-4)}^{g-1}$ (resp.\ $M_{(g-4)}^{g-2}$) is the set $\Mon_{(g-4)}^{g-1}$ (resp.\ $\Mon_{(g-4)}^{g-2}$), which is in bijection with the set of partitions of~$g-4$ with $3$ (resp.~$2$) parts.

For $g=9$ the map 
$$
\cD \colon M_{(5)}^8 = \mQ \cdot x_2^2 x_4 \oplus \mQ \cdot x_2 x_3^2 \longrightarrow M_{(5)}^7 = \mQ \cdot x_2 x_5 \oplus \mQ \cdot x_3 x_4
$$
is given by the matrix $\left({30 \atop 6}\; {20 \atop 20}\right)$; so it is surjective. Similarly, for $g=10$ the map $\cD \colon M_{(6)}^9 \to M_{(6)}^8$ is given, for the natural monomial bases of the spaces involved, by the matrix 
$$
\pmatrix{42 & 35 & 0\cr
6 & 15 & 60\cr
0 & 10 & 0 \cr}
$$
which has full rank.

Assuming now that $g \geq 11$ we have the relations
$$
\eqalign{
\cD(x_2^2 x_{g-5}) &=  (g-3)(g-4)\, x_2 x_{g-4} + 6\, x_3x_{g-5}\; ,\cr
\cD(x_2 x_3 x_{g-6}) &= {g-3\choose 3}\, x_2 x_{g-4} + {g-4\choose 2}\, x_3 x_{g-5} + 10\, x_4 x_{g-6}\; ,\cr
\cD(x_2 x_4 x_{g-7}) &= {g-3\choose 4}\, x_2x_{g-4} + {g-5\choose 2}\, x_4 x_{g-6} + 15\, x_5 x_{g-7}\; ,\cr
\cD(x_3^2 x_{g-7}) &= 2 \cdot {g-4\choose 3}\, x_3 x_{g-5} + 20\, x_5 x_{g-7}\; .\cr
}
$$
One verfies by direct calculation that these relations are linearly independent; hence $x_2 x_{g-4}$, $x_3 x_{g-5}$, $x_4 x_{g-6}$ and $x_5 x_{g-7}$ are all in the image of $\cD \colon M_{(g-4)}^{g-1} \to M_{(g-4)}^{g-2}$. But this image also contains $\cD(x_2 x_i x_{g-3-i})$ for all~$i$ with $2 \leq i \leq {g-3\over 2}$, and this element is a linear combination with positive coefficients of $x_2 x_{g-4}$, $x_i x_{g-2-i}$ and $x_{i+1} x_{g-3-i}$. Using induction on~$i$ this gives the desired surjectivity of~$\cD$.

For the second assertion of the Corollary, assume $g \geq 12$. This time we need to show that $\cD^2 \colon M_{(g-5)}^{g-3} \to M_{(g-5)}^{g-1}$ is surjective.

By direct calculation we find the relations
$$
\eqalign{
\cD^2(x_2^3 x_{g-7}) 
&= 36 \cdot {g-3\choose 4}\, x_2 x_{g-5} + 36 \cdot {g-5\choose 2}\, x_3 x_{g-6} + 180 \, x_4 x_{g-7}\, ,\cr
\cD^2(x_2^2 x_3 x_{g-8}) 
&= 40 \cdot {g-3\choose 5}\, x_2 x_{g-5} + 3g \cdot {g-5\choose 3}\, x_3 x_{g-6} + 20(g-6)(g-7) \, x_4 x_{g-7}\cr
&\qquad\qquad\qquad + 420 x_5 x_{g-8}\, ,\cr
\cD^2(x_2 x_3^2 x_{g-9}) 
&= 40 \cdot {g-3\choose 6}\, x_2 x_{g-5} + 40\cdot {g-4\choose 5}\, x_3 x_{g-6} + 40 \cdot {g-6\choose 3}\, x_4 x_{g-7}\hfill\cr
&\qquad\qquad\qquad + 40\cdot {g-7\choose 2}\, x_5 x_{g-8} + 1120\, x_6 x_{g-9}\, ,\qquad\cr
\cD^2(x_2^2 x_4 x_{g-9}) 
&= 60 \cdot {g-3\choose 6}\, x_2 x_{g-5} + 12\cdot {g-5\choose 4}\, x_3 x_{g-6} + 12 \cdot {g-5\choose 4}\, x_4 x_{g-7}\hfill\cr
&\qquad\qquad\qquad + 60\cdot {g-7\choose 2}\, x_5 x_{g-8} + 840\, x_6 x_{g-9}\, ,\qquad\cr
\cD^2(x_2 x_3 x_4 x_{g-10}) 
&= 70 \cdot {g-3\choose 7}\, x_2 x_{g-5} + 30\cdot {g-4\choose 6}\, x_3 x_{g-6} + 4g \cdot {g-6\choose 4}\, x_4 x_{g-7}\hfill\cr
&\qquad\qquad\qquad + 30\cdot {g-7\choose 3}\, x_5 x_{g-8} + 70\cdot {g-8\choose 2}\, x_6 x_{g-9} + 2520\, x_7 x_{g-10}\, ,\qquad\cr
\cD^2(x_3^3 x_{g-10}) 
&= 120 \cdot {g-4\choose 6}\, x_3 x_{g-6} + 120\cdot {g-7\choose 3}\, x_5 x_{g-8} + 3360\, x_7 x_{g-10}\, .\cr}
$$
One checks that these elements span the whole space 
$$
\mQ \cdot x_2x_{g-5} + \mQ \cdot x_3 x_{g-6} + \mQ \cdot x_4 x_{g-7} + \mQ \cdot x_5 x_{g-7} + + \mQ \cdot x_6 x_{g-9} + \mQ \cdot x_7 x_{g-10}\, ,
$$
so this subspace of $M_{(g-5)}^{g-1}$ is fully contained in the image of $\cD^2 \colon M_{(g-5)}^{g-3} \to M_{(g-5)}^{g-1}$. On the other hand, for every index~$i$ with $2 < i < {g-5\over 2}$ this image also contains the element $\cD^2(x_2^2 x_i x_{g-5-i})$, which is a linear combination with positive coefficients of the elements $x_2 x_{g-5}$, $x_3 x_{g-6}$, $x_i x_{g-3-i}$, $x_{i+1} x_{g-4-i}$ and $x_{i+2} x_{g-5-i}$. {}From this it follows by induction on~$i$ that $\cD^2$ surjects to $M_{(g-5)}^{g-1}$, as claimed. \QED

\ssection{lowgExa}
{\bf Examples.\/} In Figure~2 we give the dimensions of all~$\cR^i_{(j)}$ for some low genera. The numbering scheme is the same as in Figure~1, but note that there is no part with negative level. If in box $(i,j)$ a number $d$ appears, this means that $\dim_\mQ(\cR^i_{(j)}) = d$. The unnumbered boxes correspond to the summands $\cR^i_{(j)}$ that are zero.
\topinsert
$$
\epsfbox{lowgpicture.1}
$$
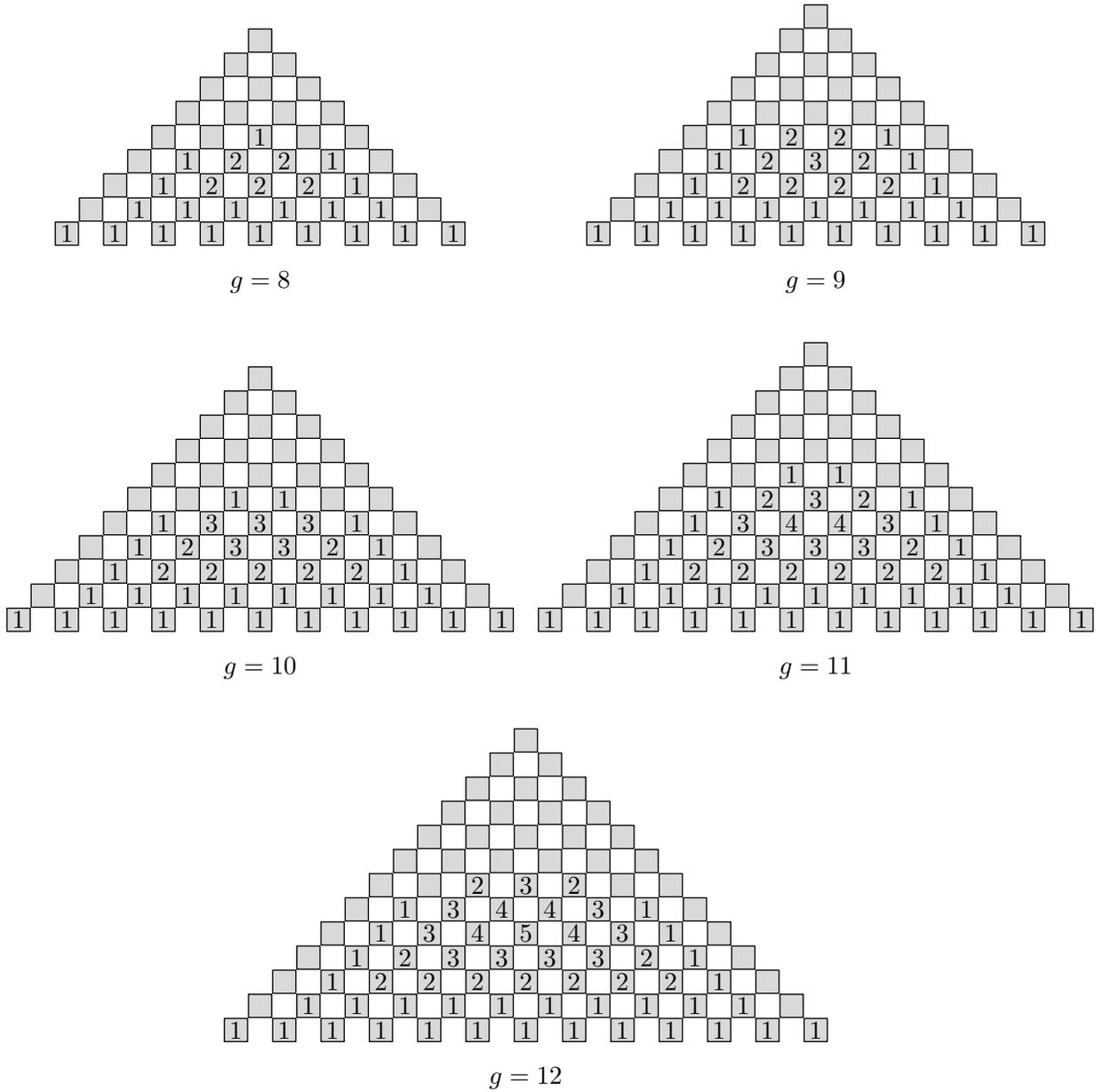
\centerline{{\bf Figure~2. The structure of the ring $\cR$ for some low values of~$g$.}}
\endinsert

Note that for $g=10$ there is a mistake in Polishchuk's list of relations in [\refn{PolUniv}], section~2.6. It is not true that $x_3x_4$ and $x_2 x_5$ are both zero, as stated there; we only have the relation $3 x_3 x_4 + 7 x_2 x_5 = 0$.

Though these tables only give the dimensions of the spaces~$\cR^i_{(j)}$, it should be clear that with little extra work one can actually write down a basis. Alternatively, using Remark~\refn{RjStructRem}(\romno2) one can give a basis for~$\cR_{(j)}$ consisting of primitive elements.

\ssection{ClosedForm}
Gerard van der Geer and Alexis Kouvidakis have informed me (personal communication of van der Geer) that they have a conjecture for the dimension of the spaces~$\cR^i_{(j)}$. Before we state their conjecture, let us introduce the notation $p_{{\bf x}}(n;{\bf y})$ for the number of partitions of~$n$ with conditions~${\bf x}$ on the number of parts, and conditions~${\bf y}$ on the parts. If we impose no condition, we omit {\bf x} or~{\bf y} from the notation. Thus, for instance, by $p_k(n;\leq l)$ we mean the number of partitions of~$n$ with $k$ parts and with all parts~$\leq l$. Conjugation of partitions interchanges the conditions ``{\bf x}'' and~``{\bf y}''; for instance, $p_k(n;\leq l) = p_{\leq l}(n;\max = k)$. Also note that $p_{\leq k}(n;\leq l) = p_k(n+k;\leq l+1)$.

The conjecture of van der Geer and Kouvidakis is that $\dim_\mQ\bigl(\cR^i_{(j)}\bigr)$ equals $p_{i-j}(i;\leq g+1-i)$, the number of partitions of~$i$ with $i-j$ parts and all parts at most $g+1-i$. We make the convention that $p_0(0;g+1) =1$. They have verified their conjecture for all~$g$ up to~$25$. Using the remarks just made, we see that the conjecture is compatible with Fourier duality.

Note that the conjecture gives that the dimension of~$\cR^i$ equals $p(i;\leq g+1-i) = p_{\leq g+1-i}(i) = p_{g+1-i}(g+1)$, and that the dimension of~$\cR$ equals $p(g+1)$, the number of partitions of $g+1$.

Van der Geer and Kouvidakis also have a strong version of the conjecture, namely that the monomials~$x^\alpha$ corresponding to the permutations of~$i$ with $i-j$ parts and all parts $\leq g+1-i$ give a basis for~$\cR^i_{(j)}$.

Though I do not know a proof for the conjecture, it is interesting to compare it with Thm.~\refn{RjStructThm}. Namely, if $2i-j \leq g$ the conjecture gives, after some rewriting, that the representation $\Sym^{g-2i+j}(\Standard)$ occurs in~$\cR_{(j)}$ with multiplicity
$$
p_{i-j}(j; \leq g-i) - p_{g+1-i}(j;\leq i-j-1)\, .
$$
The correctness of this formula suffices to prove the conjecture. Note that in our calculation of this multiplicity as the corank of the map $\cD^{g-2i+j+1} \colon M^{g-i+j+1}_{(j)} \to M^i_{(j)}$, the source space has dimension $p_{g+1-i}(j)$ and the target space has dimension $p_{i-j}(j)$.

In a first version of this paper, I hade suggested that the maps $\cD^{g-2i+j+1}$ always have the maximum possible rank. This turns out not to be true. The first counterexample occurs for $g=17$, with $i=12$ and $j=10$; in that case we look at $\cD^4 \colon M^{16}_{(10)} \to M^{12}_{(10)}$, in which source and target are both $5$-dimensional, but which has rank only~$4$.

As pointed out to me by Gerard van der Geer, the conjecture fits nicely with Brill-Noether theory. Namely, if we fix $r = i-j$ and look for the smallest~$j$ such that $\cR^{j+r}_{(j)} = 0$ then the predicted result is exactly what could be expected from geometry. To be precise, given $g$ and~$r$, let $d(g,r) = g+r-\lfloor {g\over r+1}\rfloor$ be the smallest positive integer~$d$ such that the general curve of genus~$g$ has a~$g^r_d$. By the main results in [\refn{vdGK}] and~[\refn{Herbaut}], see the next section, having a~$g^r_d$ gives a relation in $\cR^{j+r}_{(j)}$ for all $j \geq d-2r+1$. Now observe that the conjecture of van der Geer and Kouvidakis predicts that $\cR^{j+r}_{(j)} = 0$ if and only if $j \geq d(g,r) - 2r + 1$. For $r=1$ this is Polishchuk's result on the vanishing of the classes~$p_n$ for $n \geq {g\over 2} + 1$. For other small~$r$, say $r=2$ and~$3$, this can be proven with arguments as in~\refn{j=gmin4}, but as it seems difficult to get a general result with this method, we shall not give the details here.
\goodbreak


\section{Cycle relations for curves with a given linear system}{CRModAlg2}

\ssection{RelsMAIntro}
As before we consider a non-singular curve~$C$ of genus~$g$ over some algebraically closed field. The starting point in this section is the following result of Herbaut. The result was proven in~[\refn{Herbaut}] in a different formulation. Soon thereafter the result as stated here was proven by van der Geer and Kouvidakis in~[\refn{vdGK}]. It was shown by Zagier in an appendix to~[\refn{vdGK}] that the two results are actually equivalent. We shall further refine this result in Section~\refn{CRinChow}.

\ssection{HvdGK}
{\bf Theorem. {\rm (Herbaut [\refn{Herbaut}], van der Geer-Kouvidakis [\refn{vdGK}])}\/}\enspace {\it If the curve $C$ has a $g^r_d$ then
$$
\sum_{{m_1,\ldots,m_r\atop m_1+\cdots+m_r = i}} m_1!\, \cdots m_r!\cdot p_{m_1} \cdots p_{m_r} = 0  \eqlabel{HerbautRels}
$$
in $A(J)$ for all $i > d-r$.}
\Bskip

\ssection{HerbDunstable}
The idea of this section is that from Herbaut's relations we get other ones by (repeated) application of the operator~$\cD$. In general this gives many new relations, and it seems difficult to give an explicit set of generators for the ideal of relations thus obtained. So instead we shall focus on the relations we get by applying $\cD^{r-1}$ to~(\refn{HerbautRels}). What is clear a priori is that this gives a relation of the form $c \cdot p_{i+1-r} = 0$ in~$A(J)$, for some comstant~$c$ depending on $g$, $i$ and~$r$. The constant $c$ that appears turns out to be of a combinatorial nature. See Thm.~\refn{DrminHerb} and~\refn{PsiExplain} below for the precise result and some further explanation.

\ssection{BsandCs}
Given integers $i>s$, define elements $B(i,s)$ and $C(i,s)$ of~$R$ by
$$
B(i,s) := \sum_{{m_1,\ldots,m_s \geq 1\atop m_1+\cdots+m_s =i}} m_1!\, \cdots m_s!\cdot x_{m_1} \cdots x_{m_s}
\quad\hbox{and}\quad
C(i,s) := \sum_{{m_1,\ldots,m_s \geq 2\atop m_1+\cdots+m_s =i}} m_1!\, \cdots m_s!\cdot x_{m_1} \cdots x_{m_s}\, .
$$
We have
$$
B(i,s) = \sum_{a=0}^s \; {s\choose a}\; x_1^a \, C(i-a,s-a)\, .
$$
(To see this use that $B(i,s)$ is the coefficient of $t^i$ in $\bigl(\sum_{m \geq 1}\,  m!\, x_m t^m\bigr)^s$, and separate the terms with $m=1$ from those with $m \geq 2$.)

\ssection{Drhomin1C}
{\bf Lemma.\/}\enspace {\it We have
$$
\cD^{s-1}\bigl(C(i,s)\bigr) = {i! \, (i-s-1)!\over (i-2s)!} \cdot x_{i+1-s}
$$ 
for all integers $i \geq 0$ and $s \geq 1$.}
\Dskip

\noindent
{\it Proof.\/} Change to the variables $y_n := (n+1)!\ x_{n+1}$; then $\cD = \cE_1 + \cE_2$ with
$$
\cE_1 := \sum_{m,n \geq 1}\; y_{m+n} \, \partial_m \partial_n
\qquad\hbox{and}\qquad
\cE_2 := \sum_{m,n \geq 1}\; {m+n\over 2}\cdot y_{m+n}\,  \partial_m \partial_n\, .
$$
Letting
$$
C^\prime(k,s) := \sum_{{n_1,\ldots,n_s\atop n_1 + \cdots + n_s = k}}\; y_{n_1} \cdots y_{n_s}
$$
it is clear that $\cD^{s-1}\bigl(C^\prime(k,s)\bigr) = c \cdot y_k$ for some rational number~$c$. Our goal is to prove that $c = (k+s)!\, (k-1)!/ (k-s)!\, (k+1)!$. This can be checked by evaluation at the vector ${\bf 1} := (1,1,\ldots)$.

On $\mQ[y_1,y_2,\ldots]$ consider gradings~$d$ and~$\delta$ given by $d(y_i) = 1$ and $\delta(y_i) = i$. Suppose $f$ is homogeneous for both. Then one finds without difficulty (reduce to monomials) that $\cE_1(f)\bigl({\bf 1}\bigr) = d(d-1) \cdot f({\bf 1})$ and $\cE_2(f)\bigl({\bf 1}\bigr) = \delta(d-1) \cdot f({\bf 1})$. Note further that $\cE_1$ and~$\cE_2$ both decrease~$d$ by~$1$ and preserve~$\delta$. Hence, by induction on~$n \leq s-1$,
$$
\cD^n\bigl(C^\prime(k,s)\bigr)({\bf 1}) = {(k-1)!\, (k+s)!\over (k-s)!\, (k+s-n)!\, (s-1-n)!} \, ,
$$
and taking $n = s-1$ gives the result. \QED   

\ssection{FactorDef}
{\bf Definition.}\enspace Given integers $g$, $i$ and~$r$, let
$$
\Psi(g,i,r) :=  \sum_{a=0}^{r-1}\, {i-g\choose a}\, {i-a \choose r-a}\, {i-r-1\choose r-1-a}
$$
which is a polynomial of degree $r-1$ in~$g$ and of degree $2r-1$ in~$i$.

\ssection{DrminHerb}
{\bf Theorem.}\enspace {\it Suppose the curve $C$ has a $g^r_d$. Then $p_k = 0$ in~$A(J)$ for all $k > d+1-2r$ with $\Psi(g,k+r-1,r)\neq 0$. Equivalently (by Fourier duality), $[C]_{(l)} = 0$ in~$A(J)$ for all $l > d-2r$ with $\Psi(g,l+r,r)\neq 0$.}
\Dskip

\Proof
Herbaut's theorem tells us that $\phi\bigl(B(i,r)\bigr) = 0$ in~$A(J)$ for all $i > d-r$. Using Lemma~\refn{Dnp1aalpha} we find
$$
\displaylines{
\qquad\cD^{r-1}\bigl(B(i,r)\bigr) = \sum_{a=0}^r\; {r\choose a}\; \cD^{r-1}\bigl(x_1^a C(i-a,r-a)\bigr)\hfill\cr
\hfill= \sum_{a=0}^r \sum_{b=0}^{\min(a,r-1)}\; {r\choose a} {(r-1)!\over (r-1-b)!} {a!\over (a-b)!} {i-a-g+b\choose b} \cdot x_1^{a-b} \cD^{r-1-b}\bigl(C(i-a,r-a)\bigr)\, .\qquad\cr}
$$
But $\cD^{r-1-b}\bigl(C(i-a,r-a)\bigr) = 0$ if $b<a$. So
$$
\eqalign{\cD^{r-1}\bigl(B(i,r)\bigr) 
&= \sum_{a=0}^{r-1}\; {r\choose a}\; {(r-1)!\, a!\over (r-1-a)!} {i-g\choose a} {(i-a)!\, (i-r-1)!\over (i-2r+a)!} \cdot x_{i+1-r}\cr
&= r!\, (r-1)!\, (i-r)! \cdot \Psi(g,i,r) \cdot x_{i+1-r}\, .\cr}
$$
Setting $k = i+1-r$ we obtain the theorem. \QED

\ssection{CompareHerbaut}
{\bf Remark.} In~[\refn{Herbaut}], Thm.~4, Herbaut obtains a result that is very similar to our Thm.~\refn{DrminHerb}. Let us explain the relation between the results.

Herbaut gets as conclusion that $[C]_{(l)} = 0$ for all $l > d-2r$, provided that a certain combinatorical factor $A(r,d,g)$ is non-zero. In our result, the condition is that $\Psi(g,k+r-1,r)$ is non-zero; this is a condition on the triple $(g,i,k)$, and not on~$d$. The combinatorical formulas involved are related by the identity $\Psi(g,d-r+1,r) = (d-2r+2) \cdot A(r,d,g)$. So if $d-2r+2 > 0$ (the only case of interest), our result implies Herbaut's. On the other hand, even if $A(r,d,g) =0$, in which case Herbaut's theorem gives nothing, we can still apply our result to obtain the vanishing of the classes~$[C]_{(l)}$ for $l > d+1-2r$.

\ssection{PsiExplain}
At first sight the theorem seems to good to be true. For instance, if we apply it with the canonical linear system $|K_C|$ then we get relations of the form $\hbox{const} \cdot p_k = 0$ for all $k \geq 2$. However, by Ceresa's theorem in~[\refn{Ceresa}] we know that in general $p_2$ is not algebraically zero. (Note that $p_k = 0$ for $k\geq 2$ implies that $p_{k+1} =0$; see Herbaut~[\refn{Herbaut}], Thm.~5.) The point is of course that the factor $\Psi(g,k+r-1,r)$ may be zero.

It should be stressed therefore that the theorem is of interest only for special linear systems. {}From Cor.~\refn{pnnonvanish} we can deduce that $\Psi(g,i,r) = 0$ if $i \geq g$ and $r \geq {g+1\over 2} + i + 1$; this means that for non-special linear systems we get nothing new. 

\ssection{g2d}
{\bf Examples.} (a) Suppose we have a $g^2_d$. Theorem~\refn{HvdGK} gives that
$$
\sum_{{m_1,m_2\atop m_1+m_2 = i}} m_1!\, m_2!\, \cdot p_{m_1} p_{m_2} = 0
$$
in $A(J)$, for all $i \geq d-1$. Applying $\cD$ to this, we get
$$
(i-1)!\; \bigl(i(i-1) - 2g\bigr) p_{i-1} = 0\, .
$$
in $A(J)$ for all $i \geq d-1$. So if $g \neq (d-1)(d-2)/2$, which means that the $g^2_d$ is not very ample, then we find that $p_i = 0$ for all $i \geq d-2$; but if $g = (d-1)(d-2)/2$ then we only get that $p_i = 0$ for $i \geq d-1$. Note that this last conclusion is also obtained from geometry, using that if $C$ has a $g^2_d$ then also it has a $g^1_{d-1}$.

(b)\enspace
Next suppose we have a $g^3_d$. This gives that
$$
\sum_{{m_1,m_2,m_3\atop m_1+m_2+m_3 = i}} m_1!\, m_2!\, m_3!\, \cdot p_{m_1} p_{m_2} p_{m_3} \;\sim_\alg 0
$$
for all $i \geq d-2$. Applying $\cD^2$ to this, we get that $(i-2)!\, p_{i-2}$ times the factor
$$
\Psi(g,i,3) = i(i-1)^2(i-2) - 6g (i^2-3i+3) + 6g^2
$$
is algebraically equivalent to zero for all $i \geq d-2$. A computer calculation gives the following list of pairs $(g,i)$ with $i < g \leq 10000$ for which $\Psi(g,i,3) = 0$:
$$
(8,5)\; , (10,7)\; , (21,7)\; , (25,11)\; , (66,11)\; .\eqlabel{giexcept}
$$
So unless we hit one of the rare pairs $(g,i)$ for which $\Psi(g,i,3) = 0$ we get the conclusion that $p_{i-2} = 0$ for all $i \geq d-2$. If this does not work we can use that $C$ also has a $g^2_{d-1}$, and then try to use (a) to conclude that $p_k = 0$ for all $k \geq d-3$. 

Most pairs $(g,i)$ in~(\refn{giexcept}) are geometrically irrelevant. Note that we still have the condition that $d \leq i+2$. If we check which of the listed pairs satisfy Castelnuovo's bound, the only cases left are curves of genus~$10$ with a~$g^3_{10}$, and curves of genus~$25$ with a~$g^3_{13}$. I do not know if these actually occur.  

(c)\enspace
Next suppose we have a $g^4_d$. Applying $\cD^3$ to Herbaut's relations, we get that $(i-3)!\, p_{i-3}$ times the factor:
$$
\Psi(g,i,4) = i(i-1)^2(i-2)^2(i-3) - 12g (i^4 - 8 i^3 + 26 i^2 - 43 i + 34) + 36g^2 (i^2-5i+8) - 24 g^3
$$
is zero in $A(J)$, for all $i \geq d-3$. Again we look for pairs $(g,i)$ with $i < g$ for which the factor vanishes. The following is the list of such pairs with $g < 1000$:
$$
(5,1)\; , (10,6)\; , (12,8)\; , (14,10)\; , (16,9)\; , (190,38)\; .
$$
Among these, the only cases that are geometrically relevant (in the sense that they meet Castelnuovo's bound) are curves of genus~$12$ with a~$g^4_{11}$, curves of genus~$14$ with a~$g^4_{12}$ or a~$g^4_{13}$, and curves of genus~$190$ with a $g^4_d$ with $37 \leq d \leq 41$. Again I do now know if such curves actually exist.

\ssection{Noconverse}
{\bf Remark.} The general pattern we see is that the existence of special divisors on~$C$ leads to relations between the generators~$p_n$ of the tautological ring. One may ask if there is a converse to this. See for instance Herbaut~[\refn{Herbaut}], just after Thm.~4. However, it seems to us that one cannot expect such a converse, at least not in a naive way. The point is that for varieties over~$\Qbar$ one expects that the filtration $\Fil^\gdot \Chow(J)$ introduced in~\refn{picture} satisfies $\Fil^2 = 0$. In particular, for curve~$C$ over~$\Qbar$ we should have $p_n = 0$ for all $n \geq 3$. On the other hand, the general curve of genus~$g$ over~$\Qbar$ is also general in the sense of Brill-Noether theory.

As yet this is of course only speculation, but it leads to the interesting question whether one can obtain relations between the classes~$p_n$ from the assumption that there is a non-constant map $C \to \mP^1$ with at most three critical values.


\section{The big tautological ring}{BigTaut}

\ssection{WorkinChow}
In the rest of the paper we shall be interested in cycle relations in the Chow ring of~$J$ (tensored with~$\mQ$; see our conventions in~\refn{NotConv}). Note that from now on we will consider an operator~$\cD$ on~$\Chow(J)$ that lifts the operator considered before. Also we will consider elements~$p_n$ that lift those defined in~\refn{RtoChow}.

\ssection{bigtautringdef}
{\bf Definition.} We define the {\it big tautological ring $\Taut(C) \subset \Chow(J)$\/} as the smallest $\mQ$-subalgebra of $\Chow(J)$ that contains the image of $j_* \colon \Chow(C) \to \Chow(J)$ and is stable under all operations $\cdot$, $*$, $\cF$, $n^*$ and~$n_*$. Similarly, the {\it small tautological ring $\taut(C) \subset \Chow(J)$\/} is the smallest $\mQ$-subalgebra of $\Chow(J)$ that contains the classes $[C]_{(j)}$ and is stable under all operations $\cdot$, $*$, $\cF$, $n^*$ and~$n_*$.
\Bskip

\noindent
The small and big tautological rings have the same image in $A(J) := \Chow(J)/\sim_\alg$; this image is the tautological ring $\cT(C) \subset A(J)$ considered before. In general, $\Taut(C)$ is much bigger than $\taut(C)$. To make this more precise we first recall that by the results of~[\refn{Polish}], to be briefly reviewed below, $\taut(C)$ is finitely generated. By contrast, $\Taut(C)$ contains the whole $\Chow^g(J)$; hence it also contains all $\Chow^i_{(i)}$.

\ssection{an(D)def}
As in [\refn{Polish}] we define classes $p_n \in \Chow^n_{(n-1)}(J)$ and $q_n \in \Chow^n_{(n)}(J)$ by
$$
\eqalign{
p_n &:= \cF\bigl([j(C)]_{(n-1)}\bigr) = 
\hbox{degree $n$ component of $\cF\bigl[j(C)\bigr]$},\qquad 1 \leq n \leq g;\cr
q_n &:= \cF\bigl(\theta \cdot [j(C)]_{(n)}\bigr) = 
\hbox{degree $n$ component of $\cF\bigl(\theta \cdot [j(C)]\bigr)$},\qquad 0 \leq n \leq g-1\, .\cr}
$$
In particular we have $p_1 = -\theta$ and $q_0 = g$. We set $p_n := 0$ if $n \leq 0$ or $n>g$ and $q_n := 0$ if $n < 0$ or $n \geq g$.

If $y$ is an element of $\Chow(C)$ we write
$$
a_n(y) := \cF\bigl((j_* y)_{(n)}\bigr) = 
\hbox{degree $n$ component of $\cF(j_* y)$.}
$$ 
If $D$ is a divisor on~$C$ we define $a_n(D) := a_n\bigl([D]\bigr)$, which lies in $\Chow^n_{(n)}(J)$. Note that, with $K$ the canonical class of~$C$, we have $\theta \cdot \bigl[j(C)\bigr] = {1\over 2} j_* K + [0]$; see [\refn{Polish}], Section~1. Hence 
$$
q_n = \cases{
g \cdot [J] = [J] + {1\over 2} a_0(K) & if $n=0$;\cr
{1\over 2} a_n(K) & if $n \neq 0$.\cr} \eqlabel{qnanK}
$$

\ssection{ZnDdef}
Write $\ell := c_1(\cL_J)$, with $\cL_J$ as in~\refn{situation}. As in [\refn{Polish}] we define, for $n\in \mZ_{\geq 0}$ and $a\in \Chow(J)$, an operator $A_n(a) \in \End\bigl(\Chow(J)\bigr)$ by $A_n(a)\bigl(b\bigr) = a *_n b := (p_1 + p_2)_*\bigl(\ell^n \cdot p_1^* a \cdot p_2^* b)$. For $n < 0$ we set $A_n(a) := 0$.

Let $y$ be an element of $\Chow(C)$. Given integers $m\geq 0$ and $n\geq 0$ we define $Z_{m,n}(y) := m! \cdot A_n\bigl((j_* y)_{(m+n)}\bigr)$. If $m<0$ or $n<0$ we set $Z_{m,n}(y) := 0$. If $D$ is a divisor on~$C$ we write $Z_{m,n}(D) := Z_{m,n}\bigl([D]\bigr)$. 

The operators $Y_{m,n}$ considered in [\refn{Polish}] are given by
$$
Y_{m,n} := \cases{
\id + {1\over 2}\, Z_{0,0}(K) = g \cdot \id& if $(m,n) = (0,0)$;\cr
{1\over 2}\, Z_{m,n}(K) & else.\cr}
$$
Further, Polishchuk introduces operators $\tilde X_{m,n}$ and $X_{m,n}$ given by
$$
\tilde X_{m,n} := m! \cdot A_n\bigl([C]_{(m+n-2)}\bigr)
\quad\hbox{and}\quad
X_{m,n} := \tilde X_{m,n} - mn Y_{m-1,n-1}\, .
$$
Note that the $X_{m,n}$ and $\tilde X_{m,n}$ can be non-zero only if $m+n \geq 2$. If $m<0$ or $n<0$ then we set $Y_{m,n} = \tilde X_{m,n} = X_{m,n} := 0$.

The first main result of~[\refn{Polish}] is that we have the commutation relations
$$
\eqalign{
[X_{m,n},X_{m^\prime,n^\prime}] &= (m^\prime n - m n^\prime) \cdot X_{m+m^\prime-1,n+n^\prime-1}\, ,\cr
[X_{m,n},Y_{m^\prime,n^\prime}] &= (m^\prime n - m n^\prime) \cdot Y_{m+m^\prime-1,n+n^\prime-1}\, ,\cr
[Y_{m,n},Y_{m^\prime,n^\prime}] &= 0\, ,\cr}
$$
where in the second line we restrict to elements $X_{m,n}$ with $m+n \geq 2$. These identities can be interpreted as saying that a suitable Lie algebra acts on $\Chow(J)$. 

The triple $(X_{2,0}/2,X_{1,1},X_{0,2}/2)$ defines the action of $\gsl_2$ on $\Chow(J)$ as in~\refn{sl2onchow}. The operator $\cD := {1\over 2}\, X_{2,0}$ will play an important role in what follows; note that it is given by $\cD(x) = [C]_{(0)} * x$ and that it maps $\Chow^i_{(j)}(J)$ to $\Chow^{i-1}_{(j)}(J)$. 

The second main result of~[\refn{Polish}] is that $\taut(C)$ is generated by the classes $p_m$ ($1 \leq m \leq g$) and $q_m$ ($0 \leq m \leq g-1$), and that $\cD$ acts on it as the differential operator
$$
{\textstyle {1\over 2}}\cdot \sum_{m,n \geq 1}\, {m+n \choose n} p_{m+n-1} \partial_{p_m}\partial_{p_n} + \sum_{m,n \geq 1}\, {m+n-1 \choose n} q_{m+n-1} \partial_{q_m}\partial_{p_n} - \sum_{m \geq 1}\, q_{m-1} \partial_{p_m}\, .\eqlabel{Dontaut}
$$
(Polishchuk also has explicit formulas for the other operators $X_{m,n}$ and~$Y_{m,n}$; for these we refer to~[\refn{Polish}].) Note that if we say that $\cD$ acts as the operator~$(\refn{Dontaut})$, the formal meaning is that we consider the polynomial ring $\mQ\bigl[\{p_m\}_{1 \leq m \leq g},\{q_m\}_{0 \leq m \leq g-1}\bigr]$, where we now view the $p_m$ and $q_m$ as indeterminates; then the natural homomorphism $\mQ[p_\gdot,q_\gdot] \twoheadrightarrow \taut(C) \subset \Chow(J)$ intertwines the differential operator given by~$(\refn{Dontaut})$ and the operator~$\cD$.
\Bskip

\noindent
The main purpose of this section is to extend Polishchuk's results to the big tautological ring, as follows.

\ssection{ThmTaut}
{\bf Theorem.} {\it {\rm (\romno1)} The big tautological ring $\Taut(C)$ is generated, as a $\mQ$-algebra, by the classes $p_m$ and $a_m(D)$, for $m \geq 0$ and $D$ a divisor on~$C$.

{\rm (\romno2)} Let $D_1,\ldots,D_s$ be divisors on~$C$, and let $R$ be the $\mQ$-subalgebra of $\Taut(C)$ generated by the classes $p_m$, $q_m$ and $a_m(D_i)$ for $m \geq 0$ and $i \in \{1,\ldots,s\}$. Then $R$ is stable under the operations $\cdot$, $*$, $\cF$, $n^*$, $n_*$ and~$\cD$, and $\cD$ acts on it as the differential operator
\global\advance\equano by 1\definexref{DonTaut}{\equatno}{}
$$
\displaylines{
\qquad {\textstyle {1\over 2}}\cdot \sum_{m,n \geq 1}\, {m+n \choose n} p_{m+n-1} \partial_{p_m}\partial_{p_n} + \sum_{m,n \geq 1}\, {m+n-1 \choose n} q_{m+n-1} \partial_{q_m}\partial_{p_n} \hfill\cr
\hfill + \sum_{i=1}^s \sum_{m,n \geq 1}\, {m+n-1 \choose n} a_{m+n-1}(D_i) \partial_{a_m(D_i)}\partial_{p_n} - \sum_{m \geq 1}\, q_{m-1} \partial_{p_m}\, .\qquad 
\llap{{\rm (\equatno)}}\cr}
$$}

For the proof of this result we will closely follow the arguments of~[\refn{Polish}].

\ssection{PLemma2.4}
{\bf Lemma.} {\it Let $D$ be a divisor on~$C$. Then for any $n\geq 0$ we have $Z_{0,n}(D)\bigl(x\bigr) = n! \cdot a_n(D) \cdot x$.}
\Dskip

For the proof we refer to [\refn{Polish}], Lemma~2.4.

\ssection{CommRels}
{\bf Proposition.} {\it We have the commutation relations
$$
\leqalignno{
\bigl[Z_{m,n}(D),Z_{m^\prime,n^\prime}(D^\prime)\bigr] &= 0 &{\rm (a)} \cr
\bigl[X_{m,n},X_{m^\prime,n^\prime}\bigr] &= (nm^\prime - mn^\prime) \cdot X_{m+m^\prime-1,n+n^\prime-1} &{\rm (b)} \cr
\bigl[X_{m,n},Z_{m^\prime,n^\prime}(D)\bigr] &= (nm^\prime - mn^\prime) \cdot Z_{m+m^\prime-1,n+n^\prime-1}(D) &{\rm (c)} \cr
}
$$
where in\/ {\rm (c)} we restrict to elements $X_{m,n}$ with $m+n\geq 2$.}
\Dskip

\Proof~Relations (a) and~(b) were proved in~[\refn{Polish}]; see loc.\ cit., Theorem~0.1 and the remark following Lemma~2.5. For~(c) see [\refn{PolRelSymm}], Thm.~0.1. \QED 

\ssection{ThmTautPf}
{\bf Proof of Theorem {\refn{ThmTaut}}.} Let $D_1,\ldots,D_s$ and $R$ be as in (\romno2) of the theorem. Consider the polynomial ring 
$$
\tilde R := \mQ\bigl[\{p_m\}_{m \geq 1},\{q_m\}_{m\geq 0},\{a_m(D_\nu)\}_{m\geq 0, 1\leq \nu \leq s}\bigr]
$$
and let $\pi \colon \tilde R \twoheadrightarrow R$ be the natural map. Define $\tilde\cD$ to be the differential operator on~$\tilde R$ given by~(\refn{DonTaut}).

We first prove that $R$ is stable under~$\cD$ and that the action of~$\cD$ is given by~(\refn{DonTaut}), in other words, $\pi \circ \tilde\cD = \cD \circ \pi$. The first step in the proof is to show that if $x$ and~$y$ are any of the variables $p_m$, $q_m$ or $a_m(D_\nu)$ then 
$$
\pi \circ \bigl[[\tilde\cD,x],y\bigr] = \bigl[[\cD,x],y\bigr] \circ \pi\, ;\eqlabel{Dxy=Dxy}
$$
here ``$x$'' stands for ``multiplication by~$x$'' on~$\tilde R$ (resp.\ on~$R$), and likewise for~``$y$''. Note that if (\refn{Dxy=Dxy}) holds for the pair $(x,y)$ then it also holds for $(y,x)$, as follows from the Jacobi identity, using that $[x,y] = 0$. It follows from Polishchuk's results in~[\refn{Polish}] that 
$$
\pi \circ [\tilde\cD,p_m] = [\cD,p_m]  \circ \pi \eqlabel{Dpm=Dpm}
$$
for all~$m$, as $\tilde\cD$ is the sum of the operator given by~(\refn{Dontaut}) and an operator that commutes with~$p_m$. Because of the relation~(\refn{qnanK}), it only remains to verify~(\refn{Dxy=Dxy}) for elements $x=a_m(D_i)$ and $y=a_n(D_j)$. Direct calculation gives $\bigl[[\tilde\cD,a_m(D_i)],a_n(D_j)\bigr] = 0$. On the other hand, using Lemma~\refn{PLemma2.4} and Prop.~\refn{CommRels} we find that also $\bigl[[\cD,a_m(D_i)],a_n(D_j)\bigr] = 0$. Hence we have~(\refn{Dxy=Dxy}).

By induction on the degree it follows from~(\refn{Dxy=Dxy}) that for all $F \in \tilde R$ we have
$$
\pi\Bigl(\bigl[\tilde\cD,a_m(D_\nu)\bigr](F)\Bigr) = \bigl[\cD,a_m(D_\nu)\bigr]\bigl(\pi(F)\bigr)\, . \eqlabel{DamD=DamD}
$$
Note that to start the induction we need to prove this relation for $F$ the unit element of~$\tilde R$, which maps to the class $[J]$ in~$R$. Then the LHS of~(\refn{DamD=DamD}) is zero. The RHS equals $\bigl[\cD,a_m(D_\nu)\bigr]\bigl([J]\bigr) = \cD\bigl(a_m(D_\nu)\bigr) - a_m(D_\nu) \cdot \cD\bigl([J]\bigr)$, which is zero, too, because $\cD\bigl([J]\bigr) \in \Chow^{-1}(J) = 0$ and $\cD\bigl(a_m(D_\nu)\bigr) \in \Chow^{m-1}_{(m)}(J) = 0$. Hence we can start the induction, and we get the relation~(\refn{DamD=DamD}). In particular, using~(\refn{qnanK}) we find that
$$
\pi \circ [\tilde\cD,q_m] = [\cD,q_m]  \circ \pi\, . \eqlabel{Dqm=Dqm}
$$
By another induction on the degree it follows from (\refn{Dpm=Dpm}), (\refn{DamD=DamD}) and~(\refn{Dqm=Dqm}) that $\pi\bigl(\tilde\cD(F)\bigr) = \cD\bigl(\pi(F)\bigr)$ for all $F \in \tilde R$, which is what we wanted to prove.

By definition $R \subset \Chow(J)$ is stable under intersection product, and as just proven it is stable under~$\cD$. As the Fourier transform can be calculated as $\cF = \exp(e) \cdot \exp(-f) \cdot \exp(e)$ (see~[\refn{Polish}], end of Section~1, or [\refn{PolFourStab}], Lemma~1.4), $R$ is also stable under~$\cF$, and hence also under Pontryagin product. The generators of~$R$ are homogeneous (for the usual grading by codimension), so $R$ is a graded subalgebra of $\Chow(J)$, and because for $y \in \Chow^i(J)$ we have $y \in \Chow^i_{(j)}(J)$ if and only if $\cF(y) \in \Chow^{g-i+j}(J)$, we find that $R$ is stable under Beauville's decomposition, hence also under all $n^*$ and~$n_*$. This proves (\romno2) of the theorem.

Finally, let $T^\prime \subset \Taut(C)$ be the $\mQ$-subalgebra generated by all classes $p_m$, $q_m$ and~$a_m(D)$. By the results just proven, $T^\prime$ is stable under all operators $\cdot$, $*$, $\cF$, $n^*$ and~$n_*$, and as $T^\prime$ contains the image of~$j_*$ we conclude that $T^\prime = \Taut(C)$. \QED

\ssection{ThmTautCor}
{\bf Corollary.\/} {\it The big tautological ring $\Taut(C) \subset \Chow(J)$ does not depend on the choice of the base point $x_0 \in C$.\/}
\Dskip

\Proof~If $x_0^\prime \in C$ is another base point, and $j^\prime \colon C \hookrightarrow J$ is the associated embedding, we have $j^\prime = t_\delta \circ j$ with $\delta = \bigl[O_C(x_0 - x_0^\prime)] \in J$. Hence for any $y \in \Chow(C)$ we have $\cF(j^\prime_* y) = e^{\ell_\delta} \cdot \cF(j_* y)$, where $\ell_\delta \in \Chow^1_{(1)}(J)$ is the first Chern class of $(\cL_J)_{|\{y\} \times J}$. Now use the theorem, the fact that $\Taut(C)$ contains $\Chow^1_{(1)}(J)$, plus the fact that $\Taut(C) \subset \Chow(J)$ is a bigraded subring with regard to the Beauville bigrading. \QED

\ssection{ThmTautCor2}
{\bf Corollary.} {\it Let $D$ be a divisor on~$C$, and consider the $\mQ$-subalgebra $R \subset \Taut(C)$ generated by the classes $p_m$, $q_m$ and $a_m(D)$. Then for all $m \geq 1 + {g\over 2}$ the class $a_m(D)$ lies in the ideal of~$R$ generated by the classes $q_n$ with $1 \leq n < {g+1\over 2}$.\/}
\Dskip

\Proof~Write $a_m := a_m(D)$. Let $S \subset R$ be the subalgebra generated by the classes $p_m$ for $m\geq 2$, together with all $q_m$ and $a_m$. (So the only generator we exclude is~$p_1$.) Then $S$ is stable under~$\cD$ and so is the ideal $I = S\cdot q_1 + \cdots + S\cdot q_{g-1}$. The induced operator~$\cDbar$ on $S/I$ is given by
$$
{\textstyle {1\over 2}}\cdot \sum_{m,n \geq 1}\, {m+n \choose n} p_{m+n-1} \partial_{p_m}\partial_{p_n} + \sum_{m,n \geq 1}\, {m+n-1 \choose n} a_{m+n-1} \partial_{a_m}\partial_{p_n}\, .
$$

First suppose $g$ is even. Consider the element $a_1 p_2^u$ for some $u \geq 1$. By induction we find that ${\cDbar}^n(a_1 p_2^u)$, for $0 \leq n \leq u-1$, is a linear combination with coefficients in $\mZ_{\geq 0}$ of terms $a_k p_{m_1} \cdots p_{m_{u-n}}$ such that $k+m_1+ \cdots + m_{u-n} = 2u+1-n$. In particular, $\cDbar^{u-1}(a_1 p_2^u)$ is a linear combination of terms $a_k p_m$ with $k+m=u+2$, and again applying~$\cDbar$ gives $\cDbar^u(a_1 p_2^u) = c \cdot a_{u+1}$ for some positive integer~$c$. (One readily sees that $c$ is non-zero.) Now use that for all $u \geq g/2$ we have $a_1 p_2^u = 0$ in~$S/I$.

Similarly, if $g$ is odd, we start with the relation $a_1 p_2^u p_3 = 0$ for $u \geq {g-1\over 2}$, and applying $\cDbar^u$ we obtain that $a_{u+2} = 0$ in~$S/I$.

This shows that $a_m \in I\cdot R$ for all $m \geq 1 + {g\over 2}$. Finally use [\refn{Polish}], Prop.~4.2, which tells us that already in the small tautological ring $\taut(C) \subset R$ the ideal $(q_1,\ldots,q_{g-1})$ is generated by the classes~$q_n$ with $n < {g+1\over 2}$. \QED


\section{Cycle relations in the Chow ring}{CRinChow}
 
\ssection{CycRelIntro}
As before, let $C$ be a complete non-singular curve of genus $g\geq 2$. Choose a base point $x_0 \in C$, and let $\cL$ be the Poincar\'e bundle on $C \times J$, normalised such that $\cL_{|\{x_0\} \times J}$ is trivial. 

{}From now on we assume $C$ has a base-point free~$g^r_d$, say~$\cG$. Let $\Gamma$ be one of the divisors in~$\cG$, let $V \subseteq H^0(C,\Gamma)$ be the linear subspace that gives~$\cG$, and let $\gamma \colon C \to \mP(V^\vee) \cong \mP^r$ be the morphism associated to~$\cG$. If there is no risk of confusion we simple write $\mP := \mP(V)$.

Consider the incidence variety
$$
Y := \bigl\{(P,\ell) \in C \times \mP \bigm| \ell \subset H^0(C,\Gamma-P)\}\, ,
$$
and let $q_1 \colon Y \to C$ and $q_2 \colon Y \to \mP$ be the two projections. As in van der Geer-Kouvidakis~[\refn{vdGK}] we want to apply Grothendieck-Riemann-Roch to the line bundle $\cM := (q_1 \times \id)^* \cL$ on $Y \times J$ and the morphism $(q_2 \times \id) \colon Y\times J \to \mP \times J$. Our main task is to refine their calculations so as to make them work on Chow level.

\ssection{GRRProp}
{\bf Proposition.} {\it We have
$$
\eqalign{
\chern\bigl((q_2 \times \id)_* \cM \bigr) &= \cF\bigl(j_*([C] - {\textstyle {1\over 2}}K)\bigr) - \cF\bigl(j_*([C] - {\textstyle {1\over 2}}K - [\Gamma])\bigr) \cdot \exp(-h)\, ,\cr
&= \cF\bigl(j_*[\Gamma]\bigr) - \sum_{i=1}^r \cF\bigl(j_*([C] - {\textstyle {1\over 2}}K - [\Gamma])\bigr)\cdot {(-h)^i \over i!}\, ,\cr}
$$
where we identify $\Chow\bigl(\mP\times J\bigr) = \Chow(J)\bigl[h\bigr]/(h^{r+1})$ with $h$ the hyperplane class on~$\mP$.}
\Dskip

\noindent
{\it Proof.\/} Let $E$ be the rank~$r$ vector bundle on~$C$ whose fibre at a point~$P$ is $H^0(C,\Gamma-P)$. More formally, let $\pi_1$, $\pi_2 \colon C \times C \to C$ be the projections; then $E := \pi_{2,*}\bigl(\pi_1^* O_C(\Gamma) \otimes O_{C\times C}(-\Delta)\bigr)$, which by our assumption that $\cG$ is base-point free is indeed locally free of rank~$r$, and which sits in an exact sequence
$$
0 \tto E \tto V \otimes_k O_C \tto O_C(\Gamma) \tto 0\, . \eqlabel{EVOCGamma}
$$
As $Y \to C$ is the projective bundle associated to~$E$ we have $\Chow(Y) = \Chow(C)\bigl[H\bigr]/(H^r - [\Gamma] \cdot H^{r-1})$, where $H := c_1\bigl(O_Y(1)\bigr)$. If we have classes $\alpha_j \in \Chow(C)$, almost all zero, then $q_{1,*}(\sum_{j\geq 0}\, \alpha_j H^j) = \alpha_{r-1} + [\Gamma]\cdot \alpha_r$. Also note that $O_Y(1) = q_2^* O_{\mP}(1)$, so if $h := c_1\bigl(O_{\mP}(1)\bigr) \in \Chow^1(\mP)$ is the hyperplane class on~$\mP$, we have $H = q_2^*(h)$.

Let $p_1\colon Y \times J \to Y$ and $p_2 \colon Y \times J \to J$ be the projections. GRR gives
$$
\eqalign{
\chern\bigl((q_2 \times \id)_* \cM \bigr) &= (q_2 \times \id)_* \bigl(\chern(\cM) \cdot p_1^*\Todd(Y/\mP)\bigr)\cr
&= \sum_{i=0}^r p_{2,*} \Bigl(\chern(\cM) \cdot p_1^*\bigl(\Todd(Y/\mP) \cdot H^{r-i}\bigr) \Bigr) \cdot h^i \, .\cr}
$$
We calculate $p_{2,*}$ by first pushing down to $C \times J$ and then pushing down via the projection map $\pr_2 \colon C \times J \to J$. This gives
\global\advance\equano by 1\definexref{GRRcalcul}{\equatno}{}
$$
\eqalignno{
\chern\bigl((q_2 \times \id)_* \cM \bigr) &= 
\sum_{i=0}^r \pr_{2,*} \Bigl(\chern(\cL) \cdot \pr_1^*q_{1,*}\bigl(\Todd(Y/\mP) \cdot H^{r-i}\bigr) \Bigr) \cdot h^i\cr
&= \sum_{i=0}^r \cF\Bigl(j_*\, q_{1,*}\bigl(\Todd(Y/\mP) \cdot H^{r-i}\bigr)\Bigr) \cdot h^i\, .&{\rm (\equatno)}\cr}
$$ 
Next we use that $\Todd(Y/\mP) = \Todd(Y/C) \cdot q_1^*\Todd(C) \cdot q_2^* \Todd(\mP)^{-1}$. But $q_2^* \Todd(\mP) = \Todd\bigl(V \otimes_k O_Y(1)\bigr)$, so using the exact sequences (\refn{EVOCGamma}) and
$$
0 \tto O_Y \tto (q_1^* E)\bigl(1\bigr) \tto T_{Y/C} \tto 0
$$
we find that $\Todd(Y/\mP) = q_1^* \Todd(C) \cdot \Todd\bigl(q_1^* O_C(\Gamma) \otimes O_Y(1)\bigr)^{-1}$. Now 
$$
\eqalign{\Todd\bigl(q_1^* O_C(\Gamma) \otimes O_Y(1)\bigr)^{-1} 
&= \sum_{m\geq 0} {(-1)^m \over (m+1)!} \cdot \bigl([\Gamma] + H\bigr)^m\cr
&= \sum_{m\geq 0} {(-1)^m\over (m+1)!}\cdot H^m + [\Gamma] \cdot \sum_{m \geq 1}  {(-1)^m \cdot m\over (m+1)!}\cdot H^{m-1}\, ,\cr}
$$
so we find
$$
q_{1,*}\bigl(\Todd(Y/\mP) \cdot H^{r-i}\bigr) =  
\cases{
[\Gamma] & if $i=0$;\cr
{(-1)^{i-1} \over i!} \cdot \bigl([C] - {\textstyle {1\over 2}}K\bigr) \cdot \bigl([C] -[\Gamma]\bigr) & if $i>0$.\cr}
$$
Putting this back into (\refn{GRRcalcul}) we get the proposition. \QED

\ssection{ciGammaDef}
To apply Prop.~\refn{GRRProp} we need some notation. We introduce classes $c_i(\Gamma)$ that are the Chern classes of a vector bundle for which the classes $a_n(\Gamma)$ are the components of the Chern character. Concretely, choose a divisor $\Gamma$ in the linear system~$\cG$ such that $\Gamma = Q_1 + Q_2 + \cdots + Q_d$ for distinct points $Q_i \in C$. Let $L_i$ be the restriction of~$\cL$ to $\{Q_i\} \times J$, and define $W := L_1 \oplus \cdots \oplus L_d$. We then have $a_m(\Gamma) = \chern_m(W)$. Now define
$$
c_i(\Gamma) := c_i(W)\qquad\hbox{and}\qquad c_t(\Gamma) := c_t(W)\, .
$$
By construction we have
$$
\exp\Bigl(\sum_{m \geq 1} (-1)^{m-1} (m-1)!\, a_m(\Gamma) t^m \Bigr) = c_t(\Gamma)\, .
$$

\ssection{ApplyGRRProp}
Consider the situation as in~\refn{CycRelIntro}. Let $\cN$ be the vector bundle of rank~$d$ on $\mP \times J$ given by $\cN := (q_2 \times \id)_* \cM \otimes O_{\mP}(1)$. Our GRR calculation gives
$$
\chern(\cN) = \cF\bigl(j_*[\Gamma]\bigr) + 
\sum_{n\geq 1} \cF\bigl(j_*([C] - {\textstyle {1\over 2}} K)\bigr) \cdot{h^n\over n!}\, ,
$$
so
$$
\chern_0(\cN) = a_0(\Gamma) = d\quad\hbox{and}\quad 
\chern_m(\cN) = a_m(\Gamma) + \sum_{n=1}^m {p_{m-n} - {\textstyle {1\over 2}}a_{m-n}(K)\over n!} h^n\quad \hbox{for $m \geq 1$.}
$$
The non-trivial information we have is that
$$
c_t(\cN) = \exp\Bigl(\sum_{m \geq 1} (-1)^{m-1} (m-1)!\, \chern_m(\cN) t^m \Bigr)
$$
is a polynomial in~$t$ of degree $\leq d$. The RHS equals
$$
c_t(\Gamma) \cdot \exp\Bigl(\sum_{1 \leq n < m}\, {(-1)^{m-1} (m-1)!\, p_{m-n} h^n t^m\over n!} \Bigr) \cdot \exp\Bigl(\sum_{1 \leq n \leq m}\, {(-1)^{m-1} (m-1)!\, a_{m-n}(K) h^n t^m\over 2\cdot n!} \Bigr)\, .
$$
(Note that we may take $n < m$ in the second factor, as $p_0 =0$.) Separating terms according to their type in the Beauville decomposition this gives that for every $j \geq 0$ the expression
$$
c_t(\Gamma) \cdot \Bigl(\sum_{1 \leq n < m}\, {(-1)^{m-1} (m-1)!\, p_{m-n} h^n t^m\over n!} \Bigr)^j \cdot \exp\Bigl(\sum_{1 \leq n \leq m}\, {(-1)^{m-1} (m-1)!\, a_{m-n}(K) h^n t^m\over 2\cdot n!} \Bigr) \eqlabel{exprthatispol}
$$
is a polynomial in~$t$ of degree at most~$d$.

\ssection{Takej=r}
In the relation we have obtained we take $j=r$. Note that $h^{r+1} = 0$. Hence
$$
\eqalign{
\Bigl(\sum_{1\leq n < m} \, {(-1)^{m-1} (m-1)!\, p_{m-n} h^n t^m\over n!} \Bigr)^r &= \Bigl(\sum_{m\geq 2} (-1)^{m-1} (m-1)!\, p_{m-1} h t^m \Bigr)^r\cr
&= \bigl(\sum_{m\geq 1} (-1)^m m!\, p_m t^m \bigr)^r \cdot h^r t^r\, .\cr}
$$
We define $B(i,s)$ as $(-1)^i$ times the coefficient of~$t^i$ in $\bigl(\sum_{m\geq 1} (-1)^m m!\, p_m t^m \bigr)^s$; so
$$
B(i,s) = \sum_{{m_1,\ldots,m_s\atop m_1+\cdots+m_s =i}} m_1!\, \cdots m_s!\cdot p_{m_1} \cdots p_{m_s}\, .
$$
Then the relation we find is that $c_t(\Gamma) \cdot \sum_{i \geq 1} (-1)^i B(i,r) t^i$ is a polynomial in~$t$ of degree at most $d-r$. In particular, this gives the following result.

\ssection{ThmModRat} 
{\bf Theorem.} {\it Let $C$ be a curve with a $g^r_d$. Let $\Gamma$ be a divisor in the linear system, and define elements $c_n(\Gamma)$ as above. Then we have
$$
B(i,r) = c_1(\Gamma) \cdot B(i-1,r) - c_2(\Gamma)\cdot B(i-2,r) + \cdots + (-1)^i c_{i-1}(\Gamma) \cdot B(1,r)\, .
$$
for all $i > d-r$.}
\Bskip

\noindent
Note that it is no loss of generality to assume that the $g^r_d$ is complete and base-point free. Passing to the quotient of $\Chow(J)$ modulo algebraic equivalence we recover Theorem~\refn{HvdGK}.

\ssection{HowFurther}
The relations in Thm.~\refn{ThmModRat} are only the tip of the iceberg. There is an abundance of further relations. There are at least three methods we can use: (\romno1)~Take other values for~$j$ in~(\refn{exprthatispol}); (\romno2)~Apply the operator~$\cD$ to relations we have found; (\romno3)~Use that if the curve has a~$g^r_d$ then it also has a~$g^{r-1}_{d-1}$. The main difficulty is to extract manageable information, and here we have little to offer. So we confine ourselves to some simple examples of how we can get more relations.

\ssection{Takej=r-1}
{\bf Example.} We consider the relations we obtain by taking $j=r-1$ in~(\refn{exprthatispol}). We have
$$
\Bigl(\sum_{1\leq n < m} \, {(-1)^{m-1} (m-1)!\, p_{m-n} h^n t^m\over n!} \Bigr)^{r-1} = A_1 \cdot h^{r-1} t^{r-1} + A_2 \cdot h^r t^r\, ,
$$
with
$$
A_1 = \bigl(\sum_{m\geq 1} (-1)^m m!\, p_m t^m \bigr)^{r-1} = \sum_{\nu > 0} (-1)^\nu B(\nu,r-1) t^\nu\, ,
$$
and
$$
\eqalign{
A_2 &= {r-1\over 2} \cdot \bigl(\sum_{m\geq 1} (-1)^{m+1} (m+1)!\, p_m t^m \bigr) \cdot \bigl(\sum_{m\geq 1} (-1)^m m!\, p_m t^m \bigr)^{r-2}\cr
&= {r-1\over 2} \cdot \bigl(\sum_{m\geq 1} (-1)^{m+1} (m+1)!\, p_m t^m \bigr) \cdot \bigl(\sum_{\nu >0} (-1)^\nu B(\nu,r-2) t^\nu \bigr)\, .\cr}
$$
On the other hand,
$$
\exp\Bigl(\sum_{1 \leq n \leq m}\, {(-1)^{m-1} (m-1)!\, a_{m-n}(K) h^n t^m\over 2\cdot n!} \Bigr) = 1 + {\textstyle {1\over 2}} \cdot \Bigl(\sum_{m\geq 0} (-1)^m m!\, a_m(K) t^m\Bigr) \cdot h\, t \pmod{h^2}\, ,
$$
So we obtain from~\refn{ApplyGRRProp} that $c_t(W) \cdot \sum_{\nu > 0} (-1)^\nu B(\nu,r-1) t^\nu$ is a polynomial of degree $\leq d+1-r$, and that
$$
\displaylines{
\qquad c_t(W) \cdot \Bigl\{
\bigl(\sum_{\nu > 0} (-1)^\nu B(\nu,r-1) t^\nu \bigr) \cdot 
\bigl(\sum_{m\geq 0} (-1)^m m!\, a_m(K) t^m \bigr) \hfill\cr
\hfill + (r-1) \cdot \bigl(\sum_{m\geq 1} (-1)^{m+1} (m+1)!\, p_m t^m \bigr)\cdot \bigl(\sum_{\nu > 0} (-1)^\nu B(\nu,r-2) t^\nu \bigr) \Bigr\}\qquad\cr}
$$
is a polynomial of degree $\leq d-r$.

Concretely this means that 
$$
\sum_{n=0}^{i-1} (-1)^n c_n(\Gamma) B(i-n,r-1) = 0\qquad\hbox{for all $i \geq d+2-r$}
\eqlabel{B(,r-1)}
$$
and
\global\advance\equano by 1\definexref{j=rmin1second}{\equatno}{}
$$
\displaylines{
\qquad \sum_{n=0}^i c_n(\Gamma) \cdot \Bigl\{ \sum_{\nu > 0} (-1)^{n-i} (n-i-\nu)!\, B(\nu,r-1) a_{n-i-\nu}(K) \hfill\cr
\hfill + (r-1)\cdot \sum_{\nu > 0} (-1)^{n-i+1} (n-i+1-\nu)!\, B(\nu,r-2) p_{n-i-\nu}  \Bigr\} = 0 \qquad\qquad\llap{{\rm (\equatno)}}\cr}
$$
for all $i \geq d+1-r$.

\ssection{r=1}
{\bf Example.} Consider a curve with a $g^1_d$; so we take $r=1$ in the above. Thm.~\refn{ThmModRat} gives
$$
p_i = {1\over i!} \cdot \sum_{\nu=1}^{i-1}\; (-1)^{i+\nu+1}\, \nu!\, c_{i-\nu}(\Gamma)\, p_\nu\qquad \hbox{for all $i \geq d$.}
$$
Applying $\cD$ we get, using Thm.~\refn{ThmTaut},
$$
-q_{i-1} = \sum_{\nu=1}^{i-1} {(-1)^{i+\nu} \nu!\over i!} c_{i-\nu}(\Gamma) q_{\nu-1} + \sum_{\nu=1}^{i-1}\sum_{m\geq 1} {(-1)^{i+\nu+1} \nu!\over i!} {m+\nu-1\choose \nu} a_{m+\nu-1}(\Gamma) \partial_{a_m(\Gamma)}\bigl(c_{i-\nu}(\Gamma)\bigr)\, . 
$$
Now use that $\partial_{a_m(\Gamma)}\bigl(c_k(\Gamma)\bigr) = (-1)^{m-1} (m-1)!\, c_{k-m}(\Gamma)$. So we find that
$$
-q_{i-1} = \sum_{\nu=1}^{i-1} {(-1)^{i+\nu} \nu!\over i!} c_{i-\nu}(\Gamma) q_{\nu-1} + \sum_{\nu \geq 1} {(-1)^{i+1+\nu} (\nu -1)\cdot \nu!\over i!} \cdot a_\nu(\Gamma) c_{i-1-\nu}(\Gamma)
$$
for all $i \geq d$.

\ssection{r=2}
Next consider a curve with a~$g^2_d$, i.e., we take $r=2$. The first type of relation we have is the one give by Thm.~\refn{ThmModRat}, which for every $i \geq d-1$ gives an expression of $B(i,2) = \sum_{m=1}^{i-1}\; m!(i-m)!\; p_m p_{i-m}$ as a linear combination of terms $c_n(\Gamma) \cdot B(i-n,2)$.

Next we can use that the $g^2_d$ gives rise, in several ways, to a $g^1_{d-1}$. Concretely, let $\Gamma = Q_1 + \cdots + Q_d$ be one of the divisors of the~$g^2_d$, where the~$Q_i$ are distinct. For $s \in \{1,\ldots,d\}$ let $\Gamma^{(s)} := \Gamma - Q_s$. Applying Thm.~\refn{ThmModRat} to the $g^1_{d-1}$'s thus obtained we get relations
$$
p_i = {1\over i!} \cdot \sum_{\nu=1}^{i-1}\; (-1)^{i+\nu+1}\, \nu!\, c_{i-\nu}(\Gamma^{(s)})\, p_\nu \eqlabel{RelsGammai}
$$
for all $i \geq d-1$ and all $s \in \{1,\ldots,d\}$. Summing over~$s$ this gives 
$$
p_i = {1\over d \cdot i!} \cdot \sum_{\nu=1}^{i-1}\; (-1)^{i+\nu+1}\, \nu!\, (d+\nu-i) \, c_{i-\nu}(\Gamma)\, p_\nu \qquad \hbox{for all $i \geq d-1$.} \eqlabel{MoreRelspi}
$$

Thirdly we can apply what we found in~\refn{Takej=r-1}. Equation~(\refn{B(,r-1)}) gives
$$
p_i = {1\over i!} \cdot \sum_{\nu=1}^{i-1}\; (-1)^{i+\nu+1}\, \nu!\, c_{i-\nu}(\Gamma)\, p_\nu \qquad\hbox{for all $i \geq d$}
$$
which also follows from~(\refn{RelsGammai}). Equation~(\refn{j=rmin1second}) gives
$$
\sum_{n=0}^i \sum_{\nu =1}^{n-i}\; (-1)^{n-i}\, (n-i-\nu)!\, \nu!\; p_\nu\, c_n(\Gamma)\, a_{n-i-\nu}(K) = 0\, .\eqlabel{YetMoreRels}
$$
(Recall that the elements $a_l(K)$ are essentially the classes~$q_l$; see~(\refn{qnanK}).)

Finally we can apply the operator~$\cD$ to relations that we have found. For instance, applying~$\cD$ to the relation in Thm.~\refn{ThmModRat} we find a relation
$$
\eqalign{
0 = &\sum_{m=1}^{i-1}\; (-1)^{m+1}\, m\cdot (m+1)!\; c_{i-m-1}(\Gamma)\, p_m\hfill\cr
&+ 2\cdot \sum_{m=1}^{i-1}\sum_{u=0}^{i-m-1}\; (-1)^{m+u}\, u\cdot u!\, m!\; a_u(\Gamma) \, c_{i-m-u-1}(\Gamma)\, p_m\cr
&+ 2\cdot \sum_{m=1}^{i-1}\sum_{u=0}^{i-m-1}\; (-1)^{m+u}\, (u+1)!\, m!\; q_u \, c_{i-m-u-1}(\Gamma)\, p_m =0\cr}
$$
for all $i \geq d-1$. Using this identity we can express the classes $p_i$ for $i \geq d-2$ as linear combinations of the $p_{i-\nu}$ ($\nu > 0$) with coefficients in the ring $\mQ\bigl[q_m,a_m(\Gamma); m\geq 1\bigr]$. It seems that these relations (for $i \geq d-1$) do not, in general, follow from the relations we have already obtained in (\refn{RelsGammai}), (\refn{MoreRelspi}) and~(\refn{YetMoreRels}).

In conclusion, we have relations between the generators $p_m$, $q_m$ and $a_m(\Gamma)$ galore. At this stage, however, we do not have a simple set of generators for the whole ideal of relations.

\ssection{pin(q)Rem}
{\bf Remark.} A general curve of genus~$g$ has gonality $\lceil{g\over 2}\rceil + 1$. As Polishchuk shows in [\refn{Polish}], Prop.~4.2, in the small tautological ring $\taut(C)$ all classes $p_m$ with $m \geq (g/2) +1$ lie in the ideal generated by the classes~$q_n$ for $n \geq 1$. So one might guess that if the curve has a~$g^1_d$, all classes $p_n$ with $n \geq d$ should lie in the ideal $\bigl(\{q_n\}_{n \geq 1}\bigr)$, at least for a suitable choice of a base point. This, however, is probably too optimistic, except in special situations. (The main difficulty in disproving this lies in the dependence on the base point.)

As an example of a special situation, suppose the curve has a $g^r_d$ of divisors~$\Gamma$ that are (rational) multiples of the canonical class~$K$. In this case all $c_i(\Gamma)$ are in the ring $\mQ[q_l; l\geq 1]$. We find, for all $i > d-r$, that $B(i,r)$ is an element of the small tautological ring, and lies in the ideal generated by the classes~$q_l$. We can push this a bit further by the method of Section~\refn{CRModAlg2}; the conclusion we get, still assuming that $[\Gamma] \in \mQ \cdot K$, is that $p_k$ lies in the ideal generated by the~$q_l$, for all $k > d+1-2r$ with $\Psi(g,k+r-1,r) \neq 0$. A similar result has been obtained independently by Fu and Herbaut in~[\refn{FuHerbaut}].
\goodbreak


\references
\nobreak\Bskip\nobreak

\reflabel{BeauvQuelq}
A.~Beauville, {\it Quelques remarques sur la transformation de Fourier dans l'anneau de Chow d'une vari\'et\'e ab\'elienne\/}, in: Algebraic geometry (Tokyo/Kyoto, 1982), M.\ Raynaud and T.\ Shioda, eds., Lecture Notes in Math.\ {\bf 1016}, Springer, Berlin, 1983, pp.\ 238--260.

\reflabel{BeauvACJ}
A.~Beauville, {\it Algebraic cycles on Jacobian varieties\/}, Compositio Math.\ {\bf 140} (2004), 683--688.

\reflabel{Ceresa}
G.~Ceresa, {\it $C$ is not algebraically equivalent to $C\sp{-}$ in its Jacobian\/},  Ann.\ of Math.\ {\bf 117} (1983), 285--291.

\reflabel{ColvGeem}
E.~Colombo and B.\ van Geemen, {\it Note on curves in a Jacobian\/}, Compositio Math.\ {\bf 88} (1993), 333--353.

\reflabel{Fakh}
N.~Fakhruddin, {\it Algebraic cycles on generic abelian varieties\/}, Compositio Math.\ {\bf  100} (1996), 101--119. 

\reflabel{FuHerbaut}
B.~Fu and F.~Herbaut, {\it On the tautological ring of a Jacobian modulo rational equivalence\/}, preprint {\tt arXiv:0706.2814}.

\reflabel{vdGK}
G.~van der Geer and A.~Kouvidakis, {\it Cycle relations on Jacobian varieties\/}, Compositio Math.\ {\bf 143} (2007), 900--908.

\reflabel{Herbaut}
F.~Herbaut, {\it Algebraic cycles on the Jacobian of a curve with a linear system of given dimension\/}, Compositio Math.\ {\bf 143} (2007), 883--899.

\reflabel{Ikeda}
A.~Ikeda, {\it Algebraic cycles and infinitesimal invariants on Jacobian varieties\/}, J.\ Algebraic Geom.\ {\bf 12} (2003), 573--603. 

\reflabel{JannMotSh}
U.~Jannsen, {\it Motivic sheaves and filtrations on Chow groups\/}, in: Motives (Seattle, WA, 1991), U.\ Jannsen, S.\ Kleiman and J-P.\ Serre, eds., Proc.\ Symp.\ Pure Math.\ {\bf 55}, AMS, Providence, RI, 1994, pp.\ 245--302.

\reflabel{Kunne}
K.~K\"unnemann, {\it A Lefschetz decomposition for Chow motives of abelian schemes\/},
Inventiones Math.\ {\bf 113} (1993), 85--102. 

\reflabel{PolUniv}
A.~Polishchuk, {\it Universal algebraic equivalences between tautological cycles on Jacobians of curves\/},  Math.\ Z.\ {\bf 251} (2005), 875--897.

\reflabel{Polish}
A.~Polishchuk, {\it Lie symmetries of the Chow group of a Jacobian and the tautological subring\/},  J.\ Alg.\ Geom.\ {\bf 16} (2007), 459--476.

\reflabel{PolRelSymm}
A.~Polishchuk, {\it Algebraic cycles on the relative symmetric powers and on the relative Jacobian of a family of curves. I\/}, preprint {\tt arXiv:0704.2848}.

\reflabel{PolFourStab}
A.~Polishchuk, {\it  Fourier-stable subrings in the Chow rings of abelian varieties\/}, preprint {\tt arXiv:}\allowbreak{\tt 0705.0772}.

\bigskip

\noindent
{\eightrm Ben Moonen, University of Amsterdam, KdV Institute, Plantage Muidergracht~24, 1018~TV Amsterdam, The Netherlands\par}
\bye

Define $\cL := (j \times \id_J)^* \cL_J$ on $C \times J$.
\bigskip

Define $\cD \colon \Chow(J) \to \Chow(J)$ by $\cD(y) = [C]_{(0)} * y$. If $y \in \Chow^i_{(j)}(J)$ then $\cD(y) \in \Chow^{i-1}_{(j)}$. As shown by Polishchuk in~??, the operator~$\cD$ preserves the tautological subring $\cT(C) \subset \Chow(J)$, and it acts on $\cT(C)$ via the differential operator
$$
{1\over2} \sum_{m,n \geq 1}\, {m+n \choose n} p_{m+n-1} \partial_{p_m} \partial_{p_n} + \sum_{m\geq 1, n\geq 1} {m+n-1\choose n} q_{m+n-1} \partial_{q_m} \partial_{p_n} - \sum_{n\geq 1} q_{n-1} \partial_{p_n} \, .
$$
Note that $[C]_{(0)} = (-1)^{g-1} \cdot \cF(\theta)$, so we also have $\cD(y) = (-1)^{g-1} \cdot \cF(\theta) * y$.
\bigskip

\noindent
We first want to look at the special case of a curve~$C$ with a base-point free~$g^1_d$. Let $\Gamma \in \Chow^1(C)$ be the class of the~$g^1_d$, and let $\gamma \colon C \to \mP^1$ be the corresponding morphism. Let $P_1 + \cdots + P_d$ be a fibre of~$\gamma$ and define $E := \oplus_{i=1}^n L_{P_i}$ with $L_P := \cL_{|\{P\}\times J}$. Then $E$ is a vector bundle of rank~$d$ with $\chern(E) = \cF(j_* \Gamma)$.

We have $\Chow(\mP^1 \times J) = \Chow(J)\bigl[H\bigr]/(H^2)$, with $H = \bigl[\{\hbox{pt}\} \times J\bigr]$ the pull-back of the class of a point on~$\mP^1$.
\bigskip

\noindent
{\bf Lemma.} Let $V := (\gamma \times \id_J)_*\bigl(\cL \otimes \pr_C^* O_C(\Gamma)\bigr)$, which is a vector bundle of rank~$d$ on $\mP^1 \times J$. Then we have
$$
\chern(V) = \cF(j_* \Gamma) + \bigl(\cF[j(C)] - {1\over 2} \cF(j_* K)\bigr) \cdot H\, .
$$

In particular the lemma gives that $\chern_n(V) = \chern_n(E) + (p_{n-1} - \tilde q_{n-1})H$. Now we use the identity $c_t(V) = \exp\bigl(\sum_{n\geq 1}\, (-1)^{n-1} (n-1)!\, \chern_n(V) t^n \bigr)$. As $H^2=0$ we get the relation $c_t(V) = c_t(E) \cdot \bigl\{1 + H \cdot \bigl(\sum_{n\geq 0} (-1)^n n! (p_n-\tilde q_n)t^{n+1}\bigr)\bigr\}$. The condition this gives is that the RHS is a polynomial of degree $\leq d$; hence $c_t(E) \cdot \bigl(\sum_{n\geq 0} (-1)^n n!\, (p_n-\tilde q_n)t^n\bigr)$ is a polynomial of degree $\leq d-1$. Separating terms according to their type in the Beauville decomposition we obtain two relations, namely:
$$
\eqalign{
c_t(E) \cdot \bigl(\sum_{n\geq 0} (-1)^n n!\, p_n t^n\bigr) & \quad\hbox{is a polynomial of degree $\leq d-1$};\cr
c_t(E) \cdot \bigl(\sum_{n\geq 0} (-1)^n n!\, \tilde q_n t^n\bigr) & \quad\hbox{is a polynomial of degree $\leq d-1$}.\cr}
$$
\Dskip

\noindent
{\bf Proposition.}
Let $D$ be a divisor on~$C$. Then we have 
$$
\cD\bigl(\cF(n_* j_* [D]) \cdot \cF[j(C)] \bigr) =
- \cF\bigl(n_* j_* [D]\bigr) \cdot \cF\bigl(\theta \cdot [j(C)] \bigr)
- n \cdot \cF\bigl(n_* j_* [D]\bigr)
+ n \cdot \cF\bigl((n+1)_* j_* [D]\bigr) \eqno(*)
$$ 
for all $n\in \mZ$.
\Dskip

\Proof~It suffices to prove the formula for $D = 1\cdot Q$ with $Q \in C$. Let $R := j(Q) = \hbox{class of $O_C(Q-x_0)$}$. Let us write $n_J$ for the ``multiplication by~$n$'' map on~$J$. Then we have $n_* j_* [Q] = \bigl[n_J(R)\bigr] = \hbox{class of $O_C(n\cdot Q - n\cdot x_0)$}$.

If $\alpha$ and $\beta$ are elements of $\Chow(J)$ then we have $\cD\bigl(\cF(\alpha) \cdot \cF(\beta)\bigr) = (-1)^{g-1} \cdot \cF(\theta) * \bigl(\cF(\alpha) \cdot \cF(\beta)\bigr) = -\cF\bigl(\theta \cdot (\alpha * \beta)\bigr)$. Hence
$$
\cD\bigl(\cF(n_* j_* [D]) \cdot \cF[j(C)] \bigr) =
-\cF\bigl(\theta \cdot \bigl([n_J(R)] * [j(C)]\bigr)\bigr) =
-\cF\bigl([n_J(R)] * \bigl(t_{-n_J(R),*}(\theta) \cdot [j(C)]\bigr)\bigr)\, .
$$
In general, if $E$ is a divisor of degree~$0$ on~$C$ and $\epsilon \in J$ is the class of $O_C(E)$ then it follows from ?? that $t_{-\epsilon,*}(\theta) \cdot \bigl[j(C)\bigr] = \theta \cdot \bigl[j(C)\bigr] + j_*[E]$. Hence $t_{-n_J(R),*}(\theta) \cdot \bigl[j(C)\bigr] = \theta \cdot \bigl[j(C)\bigr] + j_*[n\cdot Q - n\cdot x_0] = \theta \cdot \bigl[j(C)\bigr] - n\cdot [0] + n \cdot [R]$. So we find
$$
\cD\bigl(\cF(n_* j_* [D]) \cdot \cF[j(C)] \bigr) =
-\cF\bigl[n_j(R)\bigr] \cdot \cF\bigl(\theta \cdot [j(C)]\bigr) - n \cdot \cF\bigl[n_J(R)\bigr] + n\cdot \cF\bigl[n_J(R)\bigr] \cdot \cF[R]\, .
$$
Finally, if $L_R$ is the restriction of $\cL_J$ to $\{R\} \times J$ then $\cF[R] = \exp\bigl(c_1(L_R)\bigr)$ and $\cF\bigl[n_J(R)\bigr] = \exp\bigl(c_1(L_{n_J(R)})\bigr) = \exp\bigl(n\cdot c_1(L_R)\bigr)$, so $\cF\bigl[n_J(R)\bigr] \cdot \cF[R] = \exp\bigl((n+1)\cdot c_1(L_R)\bigr) = \cF\bigl[(n+1)_J(R)\bigr] = \cF\bigl((n+1)_* j_* [Q] \bigr)$. \QED
\bigskip

\noindent
{\bf Corollary.} Let $\alpha_m := \cF\bigl((j_*[D])_{(m)} \bigr) = \hbox{degree $m$ component of $\cF\bigl(j_*[D]\bigr)$}$. Then $\cD(\alpha_m \cdot p_n) = -\alpha_m \cdot q_{n-1} + {m+n-1\choose n} \cdot \alpha_{m+n-1}$.
\Dskip

\Proof~We have $\alpha_m \in \Chow^m_{(m)}(J)$ and $\cF\bigl(n_* j_* [D]\bigr) = n^* \cF\bigl(j_*[D]\bigr)$, so $\cD\bigl(\alpha_m \cdot \cF[j(C)]\bigr)$ is the coefficient of~$n^m$ when we write $\cD\bigl(\cF(n_* j_* [D]) \cdot \cF[j(C)] \bigr)$ as a polynomial in~$n$. In the RHS of $(*)$ the coefficient of $n^m$ is
$$
\displaylines{\qquad
-\cF(\alpha_m) \cdot \cF\bigl(\theta \cdot [j(C)] \bigr) - \cF(\alpha_{m-1}) + \sum_{r=0}^g {r \choose m-1} \cF(\alpha_r)\hfill\cr
\hfill = -\cF(\alpha_m) \cdot \cF\bigl(\theta \cdot [j(C)] \bigr) + \sum_{s\geq 0} {m+s\choose m-1} \cF(\alpha_{m+s})\, .\qquad\qquad}
$$
The Corollary follows by separating degrees in the relation thus obtained. \QED
\bigskip

\noindent
{\bf Example.} (\romno1) Take $D = K$, so $\alpha_m = \tilde q_m$. We find $\cD(\tilde q_m p_n) = -\tilde q_m q_{n-1} + {m+n-1\choose n} q_{m+n-1}$. This agrees with...

(\romno2) Take $D = \Gamma$, so $\alpha_m = r_m$. We find $\cD(r_m p_n) = -r_m q_{n-1} + {m+n-1\choose n} r_{m+n-1}$. E.g., $\cD(p_1r_2) = (2-g)\cdot r_2$ and $\cD(p_2r_1) = r_2 - q_1 r_1$. 
\bigskip

\noindent
{\bf Base-point independent classes.}
\medskip

\noindent
If $\alpha$ and $\beta$ are elements of $\Chow(C)$ then the components of $\cF\bigl((-j)_*\alpha\bigr) \cdot \cF(j_*\beta)$ are independent of the chosen base point $x_0 \in C$. We use this to replace the classes $p_n$ and~$r_n$ by classes that are base-point independent. For the classes $q_n$ this is a little more subtle; it turns out that we can find base-point independent substitutes for the classes $q_n$ only for $n \geq 2$.

Define
$$
\eqalign{
a_n &:= \hbox{degree $n$ component of\ } {\textstyle {1\over 2}} \cF\bigl((-j)_* K\bigr) \cdot \cF\bigl[j(C)\bigr]\cr
&= \hbox{degree $n$ component of\ } (\tilde q_0 - \tilde q_1 + \tilde q_2 - \cdots)(p_1 + p_2 + p_3 + \cdots)\, .\cr}
$$
So 
$$
a_1 = (g-1)p_1\, ,\quad a_2 = (g-1)p_2 - p_1q_1\, ,\quad a_3 = (g-1)p_3 - p_2q_1 + p_1q_2\, ,\quad \ldots
$$

In a similar fashion, define
$$
\eqalign{
\rho_n &:= \hbox{degree $n$ component of\ } {\textstyle {1\over 2}} \cF\bigl((-j)_* K\bigr) \cdot \cF(j_*\Gamma)\cr
&= \hbox{degree $n$ component of\ } (\tilde q_0 - \tilde q_1 + \tilde q_2 - \cdots)(r_0 + r_1 + r_2 + \cdots)\, .\cr}
$$
So 
$$
\displaylines{
\hfill 
\rho_0 = d(g-1)\, ,\quad \rho_1 = (g-1)r_1 - d q_1\, ,\quad \rho_2 = (g-1)r_2 - r_1q_1 + d q_2\, ,\hfill\cr
\hfill \rho_3 = (g-1)r_3 - r_2q_1 + r_1q_2 - d q_3\, ,\quad \ldots\hfill\cr}
$$ 

For $n=2m$ even, define $b_n$ by
$$
\eqalign{
b_{2m} &:= \hbox{degree $2m$ component of\ } {\textstyle {1\over 4}} \cF\bigl((-j)_* K\bigr) \cdot \cF(j_* K)\cr
&= \hbox{degree $2m$ component of\ } (\tilde q_0 - \tilde q_1 + \tilde q_2 - \cdots)(\tilde q_0 + \tilde q_1 + \tilde q_2 + \cdots)\, .\cr}
$$
So
$$
b_0 = (g-1)^2\, ,\quad b_2 = 2(g-1)q_2 - q_1^2\, ,\quad b_4 = 2(g-1)q_4 -2q_1 q_3 + q_2^2\, , \ldots
$$
\Bskip

\noindent
{\bf Lemma.} We have 
$$
\cD(a_n) = \cases{-g(g-1) = {g\over g-1}\cdot  b_0 & if $n=0$\cr
-b_{n-1} & if $n$ is odd, $n>1$\cr
0 & if $n$ is even.\cr}
$$

Now we wish to define classes $b_{2m+1}$ for $m\geq 1$. First attempt: try $b_{2m+1} = \cD(a_{2m}b_2)$. Examples:
$$
b_3 = 6(g-1)^2 q_3 - 6(g-1) q_1 q_2 + 2q_1^3\, ,\quad
b_5 = -2q_1q_2^2 + (5-g)q_1^2q_3 + 2g(g-1) q_2q_3\, .
$$
So this does not give what we want: we cannot express $q_5$ in terms of~$b_5$. STOP: I think this is because my mathematica calculation is not ok -- continue with this at some later moment.

\bigskip

\noindent
We now further specialize to the case $g=5$ with a~$g^1_3$. By [], Prop.~4.2 the ring $\cT(C)$ is generated by the classes $p_1$, $p_2$, $p_3$, $q_1$ and~$q_2$. For instance, we have $0 = \cD^3(p_2^3) = -6 q_1^3 + 72 q_1q_2 - 108 q_3$, so $q_3 = -{1\over 18} q_1^3 + {2\over 3} q_1 q_2$. 

Our goal is to investigate whether $p_3$ lies in the ideal of $\cT(C)$ generated by $q_1$ and~$q_2$. This is equivalent to the question of whether $a_3$ lies in the ideal generated by $b_1 := q_1$ and $b_2$. We have the relation
$$
6p_3 = -c_2 p_1 + 2c_1 p_2 = -{\textstyle {1\over 2}} p_1 r_1^2 + p_1 r_2 + 2p_2r_1\, .
$$
Now I want to re-express this in terms of base-point independent classes. We have
$$
\eqalign{
p_1 &= a_1/(g-1)\cr
q_1 &= b_1 \cr
p_2 &= (a_2 + p_1 q_1)/(g-1) = a_2/(g-1) + a_1b_1/(g-1)^2\cr
q_2 &= (b_2 + q_1^2)/2(g-1) = b_2/2(g-1) + b_1^2/2(g-1)\cr
p_3 &= (a_3+p_2q_1-p_1q_2)/(g-1)\cr
&= a_3/(g-1) + a_2b_1/(g-1)^2 + a_1b_1^2/2(g-1)^3 - a_1b_2/2(g-1)^3\cr
r_1 &= (\rho_1 + dq_1)/(g-1) = \rho_1/(g-1) + db_1/(g-1)\cr
r_1^2 &= \rho_1^2/(g-1)^2 + 2d\rho_1 b_1/(g-1)^2 + d^2b_1^2/(g-1)^2\cr
r_2 &= (\rho_2+r_1q_1 - dq_2)/(g-1)\cr
&= \rho_2/(g-1) + b_1\rho_1/(g-1)^2 + db_1^2/2(g-1)^2 - db_2/2(g-1)^2\cr }
$$
Substituting these relations, and setting $d=3$ we find:
$$
\displaylines{\qquad 6(g-1)^2 a_3 + 6(g-1) a_2b_1 + 3a_1 b_1^2 - 3 a_1 b_2 =\hfill\cr
\hfill -{\textstyle {1\over 2}} a_1 \rho_1^2 + 3 a_1 b_1^2 + (g-1)a_1\rho_2 - {\textstyle {3\over 2}} a_1 b_2 + 2(g-1) a_2\rho_1 + 6(g-1) a_2b_1\, ,\qquad}
$$
so
$$
12(g-1)^2 a_3  = - a_1 \rho_1^2  + 2(g-1)a_1\rho_2 + 3 a_1 b_2 + 4(g-1) a_2\rho_1 \, .
$$
To get further, we want to apply the operator~$\cD$ to this relation. For the LHS this is easy: $\cD(12(g-1)^2 a_3) = - 12(g-1)^2 b_2 = (g-1)^2 \cdot \bigl(-24(g-1)q_2 + 12 q_1^2\bigr)$. The RHS is equal to $(g-1)^2$ times
$$
-(g-1) p_1 r_1^2 + 12 p_1 q_2 + 2(g-1) p_1 r_2 + 4(g-1) p_2 r_1 - 12 p_2 q_1\, .
$$
Except for the first term, we know how to calculate~$\cD$. We find that $\cD(\hbox{RHS})$ equals $(g-1)^2$ times
$$
\displaylines{
\qquad -(g-1) \cD(p_1 r_1^2) + 12 (2-g)q_2 + 2(g-1) (2-g)r_2 + 4(g-1) r_2 - 4(g-1) q_1 r_1 - 12 q_2 + 12 q_1^2\hfill\cr
\hfill = -(g-1) \cD(p_1 r_1^2) - 12(g-1) q_2 - 2(g-1)(g-4) r_2 - 4(g-1) q_1 r_1  + 12 q_1^2\, .  \qquad\cr}
$$
The relation we get from this is:
$$
\cD(p_1r_1^2) = 12q_2 - 2(g-4)r_2 -4 q_1 r_1\, .
$$
Now the whole point is to decide whether or not this is a ``new'' relation, but for this we need to calculate $\cD(p_1r_1^2)$ independently. (Yes, this works, see below.)

\vfill\eject

\noindent
I believe the proposition we had earlier generalizes as follows. (I'll give a version with 2 divisors, with more works similarly.)
\bigskip

\noindent
{\bf Proposition.}
Let $D_1$, $D_2$ be divisors on~$C$. Then we have 
$$
\displaylines{
\qquad\cD\bigl(\cF(n_{1,*} j_* [D_1]) \cdot \cF(n_{2,*} j_* [D_2]) \cdot \cF[j(C)] \bigr) =\hfill\cr
\hskip 3cm -\; \cF\bigl(n_{1,*} j_* [D_1]\bigr) \cdot \cF\bigl(n_{2,*} j_* [D_2]\bigr) \cdot \cF\bigl(\theta \cdot [j(C)] \bigr)\hfill\cr
\hskip 3cm -\; (n_1+n_2) \cdot \cF\bigl(n_{1,*} j_* [D_1]\bigr) \cdot \cF\bigl(n_{2,*} j_* [D_2]\bigr) \hfill\cr
\hskip 3cm +\; n_1 \cdot \cF\bigl((n_1+1)_* j_* [D_1]\bigr) \cdot \cF\bigl(n_{2,*} j_* [D_2]\bigr)\hfill\cr
\hskip 3cm +\; n_2 \cdot \cF\bigl(n_{1,*} j_* [D_1]\bigr) \cdot \cF\bigl((n_2+1)_* j_* [D_2]\bigr)\hfill\cr}
$$ 
for all $n_1, n_2 \in \mZ$.
\Dskip

\noindent
Write $\alpha_n$ and $\beta_n$ for the degree~$n$ components of $\cF(j_* [D_1])$ and $\cF(j_* [D_2])$, respectively. In the LHS of the equation the coefficient of $n_1^k n_2^l$ is $\cD\bigl(\alpha_k \cdot \beta_l \cdot \cF[j(C)] \bigr)$. In the RHS the coefficient of $n_1^k n_2^l$ is 
$$
\eqalign{
&-\alpha_k \beta_l \cF\bigl(\theta \cdot [j(C)]\bigr) - \alpha_{k-1}\beta_l - \alpha_k \beta_{l-1} + \bigl( \sum_{r=0}^g {r \choose k-1} \alpha_r \bigr) \cdot \beta_l + \alpha_k \cdot \bigl(\sum_{r=0}^g {r\choose l-1} \beta_r\bigr)\cr
&= -\alpha_k \beta_l \cF\bigl(\theta \cdot [j(C)]\bigr) + \sum_{s\geq 0} {k+s \choose k-1} \alpha_{k+s} \beta_l +  \sum_{s\geq 0} {l+s\choose l-1} \alpha_k \beta_{l+s}\, .\cr}
$$
Separating degrees we arrive at the relation
$$
\cD(\alpha_k\beta_l p_n) = -\alpha_k\beta_l q_{n-1} + {k+n-1\choose n} \alpha_{k+n-1} \beta_l + {l+n-1\choose n} \alpha_k \beta_{l+n-1}\, .
$$

\noindent
{\bf Example.} (\romno1) Take $D_1 = D_2 = {1\over 2}K$. We find 
$$
\cD(p_n \tilde q_k \tilde q_l) = - q_{n-1} \tilde q_k \tilde q_l + {k+n-1\choose n} \tilde q_{k+n-1} \tilde q_l + {l+n-1\choose n} \tilde q_k \tilde q_{l+n-1}\, .
$$

(\romno2) Take $D_1 = D_2 = K$. We find for instance $\cD(p_1 r_1^2) = (2-g)\cdot r_1^2$.

(\romno3) Take $D_1 = {1\over 2}K$ and $D_2 = \Gamma$. We find for instance $\cD(p_1 q_1 r_1) = (2-g)\cdot q_1 r_1$. According to calculations that I have done before, this means that $\cD(a_1 \rho_2)$ is base-point independent, as it should be.

\bigskip

\noindent
Now let's go back to the relation
$$
6p_3 = -c_2 p_1 + 2c_1 p_2 = -{\textstyle {1\over 2}} p_1 r_1^2 + p_1 r_2 + 2p_2r_1
$$
in the case of a $g^1_3$. Applying $\cD$ we get:
$$
-6q_2 = -{\textstyle {1\over 2}} \cdot (2-g) \cdot r_1^2 + (2-g) \cdot r_2 + 2\cdot r_2 - 2\cdot q_1 r_1\, .
$$
So the conclusion is that we can express $r_2$ as a linear combination of $q_2$, $r_1^2$ and $q_1r_1$. Putting this back in the first relation, we find an expression of $p_3$ as a linear combination of $p_1r_1^2$, $p_1q_2$, $p_1q_1r_1$ and~$p_2r_1$.

In order to arrive at conclusions it is maybe easier to consider the relation
$$
\eqalign{
12 a_3  &= (g-1)^{-2} \cdot \bigl( - a_1 \rho_1^2  + 2(g-1)a_1\rho_2 + 3 a_1 b_2 + 4(g-1) a_2\rho_1 \bigr) \cr
&= -(g-1) p_1 r_1^2 + 12 p_1 q_2 + 2(g-1) p_1 r_2 + 4(g-1) p_2 r_1 - 12 p_2 q_1\, ,\cr}
$$
which has the advantage that both sides are base-point independent. Applying~$\cD$ gives:
$$
\eqalign{
-12 b_2 &= -(g-1)(2-g) r_1^2 + 12(2-g) q_2 + 2(g-1)(2-g) r_2 + 4(g-1) r_2\cr
&\qquad\qquad - 4(g-1) q_1 r_1 - 12 q_2 + 12 q_1^2\cr
& = -(g-1)(2-g) r_1^2 + 12(1-g) q_2 - 2(g-1)(g-4) r_2 - 4(g-1) q_1 r_1  + 12 q_1^2
\cr
& = (g-1)^{-1} \cdot \bigl\{(g-2)\rho_1^2 - 3(g+2) b_2 - 2(g-1)(g-4) \rho_2\bigr\}  \cr}
$$
So $2(g-1)(g-4) \rho_2 = 9(g-2)b_2 + (g-2)\rho_1^2$, and substituting this gives
$$
12 (g-1)^2 a_3 =  \bigl(3 + 9(g-2)(g-4)^{-1}\bigr) \cdot a_1b_2 + \bigl(-1 + (g-2)(g-4)^{-1}\bigr) \cdot a_1\rho_1^2 + 4(g-1) a_2\rho_1\, .
$$

\bye